\documentclass{article}
\usepackage{array, longtable}
\usepackage{xcolor}
\usepackage{soul}
\usepackage{amsmath,amsthm,amssymb}
\usepackage{array, longtable}
\usepackage{xcolor}
\usepackage{soul}
\usepackage{booktabs}
\usepackage{makecell}
\usepackage{diagbox}
\usepackage{graphicx}
\usepackage{caption}
\usepackage{float}
\usepackage{bm}
\usepackage{arydshln}   
\usepackage[T1]{fontenc}
\usepackage{lmodern,amssymb,mathtools,microtype}
\usepackage[hidelinks]{hyperref}
\hypersetup{pdftitle={Unimodality of the q-super Catalan polynomials},
 pdfkeywords={super Catalan polynomials, unimodality, q-factorial ratios}}
\allowdisplaybreaks[2]
\numberwithin{equation}{section}
\newtheorem{theorem}{Theorem}[section]
\newtheorem{lemma}[theorem]{Lemma}
\newcommand{\Z}{\mathbb Z}
\newcommand{\qbinom}[2]{\genfrac{[}{]}{0pt}{}{#1}{#2}}

\usepackage{tikz}
\usetikzlibrary{arrows,shapes,positioning}
\usetikzlibrary{decorations.markings}
\tikzstyle arrowstyle=[scale=1]
\tikzstyle directed=[postaction={decorate,decoration={markings, mark=at position 0.75 with {\arrow[arrowstyle]{stealth}}}}]
\tikzstyle redirected=[postaction={decorate,decoration={markings, mark=at position 0.35 with {\arrow[arrowstyle]{stealth}}}}]

\numberwithin{equation}{section}
\begin{document}
\title{{Unimodality of the $q$-super Catalan polynomials\LARGE }}

\author{Chunyan Qin\thanks{ \ \ Email address:
chunyan\_qin@126.com (C. Qin)}}

\date{\small
$1.$ Faculty of Mathematical Sciences, Luoyang  Normal University , Luoyang,
471934, China
}

\maketitle

\begin{abstract}
We prove Allen's conjecture that the $q$-super Catalan polynomial
$A_{m,n}(q)=[2m]_q![2n]_q!/([m]_q![n]_q![m+n]_q!)$ is unimodal
for all nonnegative integers $m,n$. If $N=\max(m,n)\ge1000$, its
coefficients increase strictly from degree $N$ to
$\lfloor(m^2+n^2-mn)/2\rfloor$. The proof combines a positive
recurrence with uniform Fourier estimates for adjacent coefficient
differences and exact finite calculations.

\medskip
\textbf{MSC:} $05$A$20$ 


\end{abstract}

\section{Introduction}\label{sec:intro}
For a polynomial or formal power series $F$, write $[q^k]F$ for its
coefficient of $q^k$, with $[q^k]F=0$ for $k<0$. For $r\in\Z_{\ge0}$, put
\[
 [r]_q=1+q+\cdots+q^{r-1},\qquad
 [r]_q!=\prod_{j=1}^r[j]_q,\qquad [0]_q!=1.
\]
The empty sum is zero. The \emph{$q$-super Catalan polynomial} is
\begin{equation}\label{eq:def-supercat}
 A_{m,n}(q)=\frac{[2m]_q!\,[2n]_q!}
 {[m]_q!\,[n]_q!\,[m+n]_q!}
 =\sum_{k=0}^{D}a_kq^k,\qquad D=m^2+n^2-mn,
 \quad m,n\in\Z_{\ge0}.
\end{equation}
Warnaar and Zudilin proved that $A_{m,n}\in\Z_{\ge0}[q]$
\cite[Proposition~2]{WZ11}. Moreover,
\[
 A_{m,n}=A_{n,m},\qquad
 q^DA_{m,n}(q^{-1})=A_{m,n}(q),\qquad a_0=a_D=1.
\]
The second identity follows from
$[r]_{q^{-1}}!=q^{-\binom r2}[r]_q!$. Thus $a_k=a_{D-k}$, and
unimodality is equivalent to
$a_k-a_{k-1}\ge0$ for $1\le k\le\lfloor D/2\rfloor$.

At $q=1$, the quotient in~\eqref{eq:def-supercat} is the super
Catalan integer $(2m)!(2n)!/(m!\,n!\,(m+n)!)$; see
Gessel~\cite[\S6]{G92}. Some specializations are Gaussian polynomials.
With
\[
 \qbinom{L}{K}_q=\frac{[L]_q!}{[K]_q!\,[L-K]_q!}
 \qquad(0\le K\le L,\ K,L\in\Z),
\]
cancellation gives
\[
 A_{0,n}(q)=A_{n,n}(q)=\qbinom{2n}{n}_q,
 \qquad A_{m,2m}(q)=\qbinom{4m}{m}_q.
\]
The coefficient of $q^k$ in $\qbinom{L}{K}_q$ counts partitions of
$k$ contained in a $K$ by $L-K$ rectangle. Their unimodality has a
constructive proof due to O'Hara~\cite{OH90}. Pak and
Panova~\cite[Theorem~1]{PP13} proved strict increase from degree $1$
to the middle when both side lengths are at least $8$.
Asymptotic estimates for consecutive differences in rectangles of
bounded aspect ratio were obtained by Melczer, Panova and
Pemantle~\cite{MPP20}.

The conjecture appeared in Allen's 2009 master's thesis, as
recorded by Allen and Gheorghiciuc~\cite[p.~2 and
reference~1]{AG14}. Their weighted Dyck-path interpretation treats
$A_{2,n}(q)/(1+q^n)$~\cite[Theorem~2]{AG14}.
For $m=1$, the identity
\[
 A_{1,n}(q)=(1+q)C_n(q),\qquad
 C_n(q)=\frac{1}{[n+1]_q}\qbinom{2n}{n}_q
\]
relates the question to parity-unimodality of rational $q$-Catalan
polynomials. Xin and Zhong~\cite{XZ20} formulated that property and
proved it when the smaller parameter is at most $5$.
Stucky~\cite[Theorem~1.2]{S21} proved the general statement as a
special case of his theorem on rational $q$-Schr\"oder polynomials.
In particular, if $C_n(q)=\sum_jc_jq^j$, its even and odd coefficient
sequences are symmetric and unimodal. Since $\deg C_n=n(n-1)$,
\[
 [q^k](1-q)A_{1,n}(q)=c_k-c_{k-2}\ge0
 \qquad\left(1\le k\le\frac{n(n-1)}2\right),
\]
which proves this specialization.

Rodelet-Causse and Tevlin~\cite{RT25} constructed Narayana
refinements with nonnegative polynomial coefficients. Tevlin~\cite{T26}
proved further coefficientwise inequalities with respect to the
refinement index. For the coefficients in $q$, Warnaar and
Zudilin~\cite[Conjecture~5]{WZ25} proposed unimodality after
multiplication by $1+q$ for a broader class of factorial ratios.
For their family~(6c), which is~\eqref{eq:def-supercat}, they stated
that the factor $1+q$ should be unnecessary
\cite[remark following Lemma~7]{WZ25}.

Our main result proves Allen's conjecture and gives a uniform range
of strict inequalities.
\begin{theorem}\label{thm:unimodality}
Let $m,n\in\Z_{\ge0}$, $N=\max(m,n)$, and $D=m^2+n^2-mn$. Then
\[
 a_k\ge a_{k-1}\qquad(1\le k\le\lfloor D/2\rfloor).
\]
If $N\ge1000$, then
\[
 a_k>a_{k-1}\qquad(N<k\le\lfloor D/2\rfloor).
\]
Consequently $A_{m,n}(q)$ is symmetric and unimodal for every
pair of nonnegative integers $m,n$.
\end{theorem}

The theorem also proves \cite[Conjecture~5]{WZ25} for family~(6c).
Indeed, after setting $a_j=0$ outside $[0,D]$,
\[
 [q^k](1-q)(1+q)A_{m,n}(q)
 =(a_k-a_{k-1})+(a_{k-1}-a_{k-2})\ge0
 \quad\left(1\le k\le\left\lfloor\frac{D+1}{2}\right\rfloor\right).
\]
When $D$ is odd, the first difference at $k=(D+1)/2$ is zero by
reciprocity.

We study the coefficient of $q^k$ in $B(q)=(1-q)A_{m,n}(q)$.
Near the middle, we extract the two adjacent coefficients on the
same circle, choosing its radius so that the logarithmic derivative
of $A_{m,n}$ equals $k-1/2$. For smaller degrees we use the
logarithmic derivative of $B$ itself. A positive recurrence controls
the Fourier integral away from zero; Euler-factor estimates control
the local integral. The two ranges overlap. Cancellation proves the
remaining inequalities in degrees at most $N$, and an exact integer
calculation treats $N\le1000$.

Section~\ref{sec:prelim} fixes the notation and records the algebraic
and analytic identities used throughout. Sections~\ref{sec:small}
and~\ref{sec:large} prove positivity at the two kinds of saddles.
Section~\ref{sec:proof} verifies coverage of every lower-half
coefficient and completes the proof.

\section{Preliminaries}\label{sec:prelim}
Throughout, $m,n,N$ are nonnegative integers. By symmetry it is
enough to consider $0\le m\le n=N$. Write
$D=m^2+N^2-mN$ and $B(q)=(1-q)A_{m,N}(q)$. For a positive real parameter $\lambda$,
the circle of integration has radius $q=e^{-\lambda}$, and we put
\[
 M=\lambda^{-1},\qquad u=N\lambda=N/M.
\]
We use $\theta$ for the coefficient-extraction angle and put
$t=N\theta$ in the compact-arc estimates. Every terminating decimal denotes
the corresponding exact rational number.

We use the positivity theorem and recurrence of Warnaar and
Zudilin~\cite[Proposition~2 and (11)]{WZ11}, together with Gaussian
unimodality~\cite{OH90}. Set
\[
 (q;q)_r=\prod_{j=1}^r(1-q^j),\qquad (q;q)_0=1,
\]
and take $\qbinom{L}{K}_q=0$ unless $0\le K\le L$.
The $q$-Pascal identity and the finite $q$-binomial identity are
\[
 \qbinom{L}{K}_q=q^K\qbinom{L-1}{K}_q+\qbinom{L-1}{K-1}_q,
 \qquad
 \prod_{j=1}^{L}(1+zq^j)
 =\sum_{K=0}^{L}z^Kq^{K(K+1)/2}\qbinom{L}{K}_q.
\]

For integers $a,p\ge0$ set $b=a+p$. In the inner sum of the cited
recurrence write its index as $h+j$. Since
$h(a+h)+(h+j)(a+h+j)=2h(a+h)+j(a+2h+j)$, we obtain
\begin{equation}\label{eq:wz-mixture}
 A_{a,b}(q)=\sum_{h=0}^{\lfloor p/2\rfloor}
 A_{a,h}(q)\qbinom{p}{2h}_q\,q^{2h(a+h)}H(b,p-2h;q),
\end{equation}
where $H(b,r;q)=\sum_{j=0}^{r}q^{j(b-r+j)}\qbinom{r}{j}_q$ for $0\le r\le b$.
For each displayed factor $F$ in a summand of \eqref{eq:wz-mixture},
nonnegativity gives $|F(qe^{i\theta})|\le F(q)$ when $q>0$.

For $1\le r\le b$, apply the $q$-Pascal identity inside $H$:
the term with $q^j\qbinom{r-1}{j}_q$ gives $H(b,r-1;q)$, and the term
with $\qbinom{r-1}{j-1}_q$, on putting $j=\ell+1$, has exponent
$b-r+1+\ell(b-r+2+\ell)$. Thus
\[
 H(b,r;q)=H(b,r-1;q)+q^{b-r+1}H(b+1,r-1;q).
\]
For $0<q\le1$, replacing $b$ by $b+1$ in $H(b,r-1;q)$
multiplies its $j$th summand by $q^j\in[q^{r-1},1]$. Hence
\begin{equation}\label{eq:H-ratio}
 (1+q^b)H(b,r-1;q)\le H(b,r;q)\le 2H(b,r-1;q).
\end{equation}
Divide the $h$th summand in \eqref{eq:wz-mixture} by $A_{a,b}(q)$
and denote the result by $W_h(q)$. Then $W_h(q)>0$ and $\sum_hW_h(q)=1$.
For any function $f$ on this finite index set, write
\[
 \mathbb E_q f(h)=\sum_{h=0}^{\lfloor p/2\rfloor}W_h(q)f(h).
\]

Writing $[j]=[j]_q$ temporarily, the two factorial ratios at
$0\le h<\lfloor p/2\rfloor$ are
\[
 \frac{A_{a,h+1}}{A_{a,h}}=\frac{[2h+2][2h+1]}{[h+1][a+h+1]},
 \qquad \frac{\qbinom{p}{2h+2}_q}{\qbinom{p}{2h}_q}
 =\frac{[p-2h][p-2h-1]}{[2h+2][2h+1]} .
\]
Their product gives
\begin{equation}\label{eq:tilted-weight-ratio}
 \frac{W_{h+1}(q)}{W_h(q)}=
 q^{2a+4h+2}
 \frac{[p-2h]_q[p-2h-1]_q}{[h+1]_q[a+h+1]_q}
 \frac{H(b,p-2h-2;q)}{H(b,p-2h;q)}.
\end{equation}
At $q=1$, $H(b,r;1)=2^r$, and cancellation of the factorials gives
\[
 W_h(1)=\binom{2b}{p}^{-1}
        \frac{b!\,2^{p-2h}}{h!\,(p-2h)!\,(a+h)!}.
\]
This is the probability of selecting exactly $h$ complete pairs
in a uniform $p$-element subset of $b$ disjoint pairs:
there are $h$ complete, $p-2h$ single and $a+h$ empty pairs,
with two choices in each single pair.

By symmetry take $0\le m\le n=N$. For $N=0$ we have $A_{0,0}=1$.
Assume henceforth that $N\ge1$. The exponent of
$1-q^j$ in $A=A_{m,N}$ is
\[
 e_j=\mathbf1_{j\le2m}+\mathbf1_{j\le2N}
      -\mathbf1_{j\le m}-\mathbf1_{j\le N}-\mathbf1_{j\le m+N},
 \qquad j\ge1.
\]
The powers of $1-q$ from the $q$-integers cancel, since
$2m+2N=m+N+(m+N)$. Thus
$A(q)=\prod_{a\in\mathcal U}(1-q^a)/\prod_{b\in\mathcal V}(1-q^b)$, where
$\mathcal U=\{j:e_j=1\}$ and $\mathcal V=\{j:e_j=-1\}$ are as follows:
\[
 \begin{array}{c|c|c}
 &\mathcal U&\mathcal V\\\hline
 N\le2m &[N+1,2m]\cup[m+N+1,2N]&[1,m]\\
 N\ge2m &[m+N+1,2N]&[1,m]\cup[2m+1,N].
 \end{array}
\]
Here $[c,d]$ in an index set means $\{j\in\Z:c\le j\le d\}$.
In the first row both sets have cardinality $m$, in the second
$N-m$. Thus $|\mathcal U|=|\mathcal V|$, $1\in\mathcal V$, and $b\le N<a$ in every pair
$(a,b)\in \mathcal U\times\mathcal V$.

Put $B(q)=(1-q)A(q)$. For $C=A$ or $B$ and real $\lambda>0$, write
\[
 L_C(\lambda)=\log C(e^{-\lambda}),\qquad
 \mu_C=-L_C',\qquad V_C=L_C'',\qquad M=1/\lambda .
\]
An equation $\mu_C(\lambda)=s$ specifies a saddle at mean parameter
$s$. The logarithm exists since $A(q),B(q)>0$ for $0<q<1$.
For $A$ these definitions extend to $\lambda=0$. In fact, if
$p_j=a_je^{-\lambda j}/A(e^{-\lambda})$, then differentiating gives
\[
 \mu_A=\sum_j jp_j,\qquad p_j'=(\mu_A-j)p_j,\qquad
 -\mu_A'=V_A=\sum_j(j-\mu_A)^2p_j>0 ;
\]
the strict inequality follows from $D>0$ and $a_0=a_D=1$.
Since $p_j=p_{D-j}$ at $\lambda=0$ and $p_0\to1$ at $+\infty$,
$\mu_A$ decreases from $D/2$ to $0$. For $B$ differentiation of
$L_B=L_A+\log(1-e^{-\lambda})$ gives
\[
 \mu_B=\mu_A-\frac1{e^\lambda-1},\qquad
 V_B=V_A-\frac{e^{-\lambda}}{(1-e^{-\lambda})^2}.
\]
In particular, positivity of $V_B$ must be established on each
interval where monotonicity of $\mu_B$ is used.

From the canceled product we have
\[
 \mu_A=\sum_{b\in\mathcal V}\frac{b}{e^{\lambda b}-1}
       -\sum_{a\in\mathcal U}\frac{a}{e^{\lambda a}-1}
 \le\sum_{j\ge1}\frac{j}{e^{\lambda j}-1}
 \le\frac{\pi^2}{6\lambda^2}.
\]
Indeed, $x/(e^x-1)$ decreases on $x>0$ because
$((1-x)e^x)'=-xe^x<0$. The sum is bounded by its integral,
and expanding $(e^x-1)^{-1}=\sum_{\ell\ge1}e^{-\ell x}$ gives
$\int_0^\infty x/(e^x-1)\,dx=\zeta(2)=\pi^2/6$,
where $\zeta(s)=\sum_{\ell\ge1}\ell^{-s}$ for real $s>1$.

For $d=2,3$ put $f_d(x)=x^d\sum_{\ell\ge1}\ell^{d-1}e^{-\ell x}$.
Summing the geometric series and its derivatives gives
\[
 f_2(x)=\left(\frac{x}{2\sinh(x/2)}\right)^2,\qquad
 f_3(x)=f_2(x)x\coth(x/2),\qquad f_2(0)=1,\quad f_3(0)=2
\]
by continuous extension. Both functions decrease: one has
$(\log f_2)'=2/x-\coth(x/2)<0$, since
$y\cosh y-\sinh y=\int_0^y t\sinh t\,dt>0$ for $y>0$.
For $y=x/2$, $(\log f_3)'\le0$ is equivalent to
$g(y)=y(3+2\sinh^2 y)-3\sinh y\cosh y\ge0$;
indeed $g(0)=0$ and $g'(y)=4\sinh y\,(y\cosh y-\sinh y)\ge0$.
For any bijection pairing $\mathcal U$ and $\mathcal V$, termwise differentiation of
the Euler product therefore gives
\begin{equation}\label{eq:paired-cumulants}
 \begin{aligned}
 V_A(\lambda)&=\lambda^{-2}\sum_{(a,b)}
       \{f_2(\lambda b)-f_2(\lambda a)\},\\
 -V_A'(\lambda)&=\lambda^{-3}\sum_{(a,b)}
       \{f_3(\lambda b)-f_3(\lambda a)\}\ge0 .
 \end{aligned}
\end{equation}
At zero tilt, $\log(\sinh y/y)=y^2/6+O(y^4)$ and
$1-e^{-js}=js\,e^{-js/2}\sinh(js/2)/(js/2)$. Put
$S(t)=\sum_{j=1}^t j^2=t(t+1)(2t+1)/6$. The factors $s$ cancel in
the quotient, so the quadratic terms in its logarithm give
\[
 V_A(0)=\frac{1}{12}\sum_{j\ge1}e_jj^2
 =\frac{S(2m)+S(2N)-S(m)-S(N)-S(m+N)}{12} .
\]
Substitution and \eqref{eq:paired-cumulants} give
\begin{equation}\label{eq:variance-origin}
 V_A(\lambda)\le V_A(0)
 =\frac{2N^3+N^2-m(N-m)(2m+N+1)}{12}
 \le\frac{2N^3+N^2}{12}\qquad(\lambda\ge0).
\end{equation}
Two normalizations will recur in the Fourier integrals:
\[
 I_G=\int_0^\infty e^{-V\theta^2/2}\,d\theta
     =\sqrt{\frac{\pi}{2V}}\quad(V>0),\qquad
 I_*(M)=\sqrt{\frac{3}{2\pi}}\,M^{-3/2}.
\]
For $a>0$ the half-Gaussian moments satisfy
\[
 J_j(a):=\int_0^\infty \theta^{2j}e^{-aV\theta^2}\,d\theta
 =\frac{2j-1}{2aV}J_{j-1}(a)\qquad(j\ge1)
\]
by integration by parts. Since $J_0(a)=I_G/\sqrt{2a}$, it follows that
\begin{equation}\label{eq:half-moments}
 \frac{J_j(a)}{I_G}=\frac{(2j-1)!!}{(2aV)^j\sqrt{2a}},\qquad j\ge1.
\end{equation}
For the complex derivatives we use the classical product
\cite[Eq.~4.36.1]{DLMF}
\begin{equation}\label{eq:sinh-product}
 \frac{\sinh(w/2)}{w/2}
 =\prod_{\ell\ge1}\left(1+\frac{w^2}{4\pi^2\ell^2}\right)
 \qquad(w\in\mathbb C),
\end{equation}
with value $1$ at $w=0$.

For $d=2,3$ the positive series gives
\[
 \int_0^\infty f_d(x)\,dx
  =\sum_{\ell\ge1}\ell^{d-1}\frac{d!}{\ell^{d+1}}
  =d!\zeta(2),\qquad \int_0^\infty f_2(x)\,dx=\frac{\pi^2}{3}.
\]
Also $f_3(x)\ge2f_2(x)$ since $(x/2)\coth(x/2)\ge1$.

Dropping numerator contributions and using the right-endpoint sum for
the decreasing function $f_2$ yields a bound valid at every real tilt:
\[
 V_B(\lambda)\le V_A(\lambda)
 \le M^2\sum_{j\ge1} f_2(j/M)\le \frac{\pi^2}{3}M^3 .
\]
Let $b_k=[q^k]B=a_k-a_{k-1}$. Coefficient extraction on
$|z|=q=e^{-\lambda}$ gives the common starting point
\[
 b_k=\frac{q^{-k}}{\pi}\int_0^\pi
 \Re\!\left[e^{-ik\theta}(1-qe^{i\theta})A(qe^{i\theta})\right]d\theta.
\]
Indeed, integrate the finite expansion of $B(qe^{i\theta})$ on
$[-\pi,\pi]$ against $e^{-ik\theta}$ and pair the summands at
$\theta$ and $-\theta$. We use either $\mu_B(\lambda)=k$ or
$\mu_A(\lambda)=k-1/2$. In each case we compare the real integral near
$\theta=0$ with the absolute integral on its complement.
The degree ranges corresponding to these choices are compared in
Section~\ref{sec:proof}.

The finite coefficient calculation in Section~\ref{sec:proof} uses
C++17 and GMP integers. The harmonic boundary calculation in
Section~\ref{sec:small} uses Python with integer interval endpoints
in $10^{-36}\Z$ and outward rounding. Exponentials and cosines are enclosed by Taylor
polynomials with explicit remainders; the inequalities that enter
the proof, including their full parameter ranges, are stated below.

\section{The saddle for adjacent coefficients}\label{sec:small}
Here $A=A_{m,N}$ and $\mu_A(\lambda)=k-1/2$. We prove
$a_k>a_{k-1}$ when $N\ge1000$ and $0<N\lambda\le2$.

\subsection{Variance and separated Euler factors}

For $\lambda>0$ and $0\le\theta\le\pi/N$, the normalized modulus of
an Euler factor is nonincreasing in its index. In $A_{m,N}$ we can
pair $\lfloor N/4\rfloor$ numerator and denominator indices so that
the two indices in each pair differ by at least $N$.
\begin{lemma}\label{lem:first-lobe-separated-pairs}
Let $0\le m\le N$, $N\ge1$, $q=e^{-\lambda}$ with $\lambda>0$,
and $0\le\theta\le\pi/N$. Put
\[
 r=\lfloor N/4\rfloor,\qquad
 a=\sin^2(N\theta/4),\qquad b=\sinh^2(N\lambda/4).
\]
Then
\[
 \left|\frac{A_{m,N}(qe^{i\theta})}{A_{m,N}(q)}\right|
 \le
 \left(1-\frac{8a(3+2b-2a)}{(3+4b)^2}\right)^{r/2}
 \le
 \exp\left\{-\frac{4ra(3+2b-2a)}{(3+4b)^2}\right\}.
\]
The same inequalities hold for $\lambda=0$ by continuity.
\end{lemma}
\begin{proof}
Let $\mathcal U,\mathcal V$ be the surviving numerator and denominator indices, and set
$s=\lfloor N/2\rfloor$, $h=\lceil 3N/2\rceil$. We count the indices
in $\mathcal V\cap[1,s]$ and $\mathcal U\cap[h,2N]$. For $m\le N/2$ the two counts are
\[
 m+\max(0,s-2m)=\max(m,s-m)\ge\lceil s/2\rceil,\qquad
 \min(N-m,2N-h+1)\ge\lceil N/2\rceil .
\]
For $m\ge N/2$ the first count is $s$. The second is
\[
 N-m+\max(0,2m-h+1)
 \ge N-\left\lfloor\frac{h-1}{2}\right\rfloor
 =\lfloor N/4\rfloor+1 .
\]
The last identity holds in each of the four residue classes modulo $4$.
Thus we may pair $r$ low denominator indices with $r$ high numerator
indices and complete the pairing arbitrarily. Denote each
numerator--denominator pair by $(j,\ell)$.

For real $x>0$ put
\[
 f(x)=\frac12\log\left(1+
       \frac{\sin^2(x\theta/2)}{\sinh^2(x\lambda/2)}\right).
\]
When $\theta>0$, the ratio
$\sin(x\theta/2)/\sinh(x\lambda/2)$ is positive and decreasing for
$0<x<2\pi/\theta$, since its logarithmic derivative satisfies
\[
 \frac\theta2\cot(x\theta/2)
 -\frac\lambda2\coth(x\lambda/2)<\frac1x-\frac1x=0.
\]
Here $z\cot z<1$ for $0<z<\pi$ and $z\coth z>1$ for $z>0$.
Thus $f$ is decreasing on $0<x\le2N$, with the boundary values filled
in continuously. When $\theta=0$, both sides equal one.

Each of the selected pairs contributes at most
$f(3N/2)-f(N/2)$ to its logarithmic modulus, and every remaining
pair contributes at most zero because $j>N\ge\ell$. Thus
\[
 \log\left|\frac{A_{m,N}(qe^{i\theta})}{A_{m,N}(q)}\right|
 =\sum_{(j,\ell)}(f(j)-f(\ell))
 \le r\bigl(f(3N/2)-f(N/2)\bigr).
\]
The triple-angle identities for sine and hyperbolic sine imply
\[
 \exp\bigl(2f(3N/2)-2f(N/2)\bigr)
 =\frac{b(3+4b)^2+a(3-4a)^2}{(a+b)(3+4b)^2}
 =1-\frac{8a(3+2b-2a)}{(3+4b)^2}.
\]
The ratio is positive for $b>0$ and at most one for $0\le a\le1/2$.
Raising it to the power $r/2$ gives the first bound. The exponential
bound follows from $1-z\le e^{-z}$. As $A_{m,N}(1)>0$, both bounds
extend to $\lambda=0$ by continuity.
\end{proof}

For the variance, differentiation turns the five factorials of
$A_{m,N}$ into five finite sums. Their continuous limits involve
$f_2(x)=(x/(2\sinh(x/2)))^2$, with $f_2(0)=1$, and its integral
$F_2(t)=\int_0^t f_2(x)\,dx$. A lower bound at $N\lambda=2$ will also
control the local Gaussian integrals at smaller tilts.
\begin{lemma}\label{lem:u2-continuum-variance}
For $0\le s\le u$ and $u\ge2$,
\[
 F_2(s)+F_2(u)+F_2(s+u)-F_2(2s)-F_2(2u)>0.57.
\]
Consequently, for every $M=\lambda^{-1}>0$, $0\le m\le N$,
and $N\lambda\ge2$, one has
\[
 \frac{V_A(\lambda)}{M^3}>\frac{57}{100}-\frac3M,
 \qquad \frac{V_B(\lambda)}{M^3}>\frac{57}{100}-\frac4M,
 \qquad B(q)=(1-q)A_{m,N}(q).
\]
\end{lemma}
\begin{proof}
Write $G(s,u)$ for the five-term expression. Since $f_2$ decreases and
$s+u\le2u$,
\[
 \frac{\partial G}{\partial u}
 =f_2(u)+f_2(s+u)-2f_2(2u)\ge0.
\]
For $0\le s\le2$ it therefore suffices to take $u=2$.
For $s\ge2$, the least permitted $u$ is $s$, and
$G(s,s)=2F_2(s)-F_2(2s)$ has derivative $2(f_2(s)-f_2(2s))\ge0$.
It remains to bound $G(s,2)$ for $0\le s\le2$.

To bound $G(s,2)$, put
\[
 \ell(x)=1-\frac{x^2}{12}+\frac{x^4}{240}-\frac{x^6}{6048}
 +\frac{x^8}{172800}-\frac{x^{10}}{5322240}
 +\frac{691x^{12}}{118879488000}-\frac{x^{14}}{5748019200},
\]
and define
\[
 L(x)=\int_0^x\ell(t)\,dt,\qquad
 U(x)=L(x)+\frac{3617x^{17}}{12093372555264000}.
\]
The partial-fraction expansion of $\coth(x/2)$, or differentiation
of the product~\eqref{eq:sinh-product}, gives
\[
 f_2(x)=1+2\sum_{j\ge1}(-1)^j(2j-1)\zeta(2j)
                    \left(\frac{x}{2\pi}\right)^{2j}
 \qquad(|x|<2\pi).
\]
For $0\le x\le4$, the ratio of the second nonconstant term to the
first is $x^2/20\le4/5$. For every subsequent ratio,
\[
 \frac{x^2}{4\pi^2}\frac{2j+1}{2j-1}
       \frac{\zeta(2j+2)}{\zeta(2j)}
 <\frac{16}{36}\frac53=\frac{20}{27}<1\qquad(j\ge2).
\]
Thus the alternating-series remainder has the sign of its first
omitted term, and integration yields
$L(x)\le F_2(x)\le U(x)$ on $[0,4]$.
Consequently, with
\[
 \mathcal R(s)=L(s)+L(2)+L(s+2)-U(2s)-U(4)-\frac{57}{100},
\]
one has $G(s,2)-57/100\ge \mathcal R(s)$ for $0\le s\le2$.
Write $\mathcal R(1+y)=\sum_{j=0}^{17}c_jy^j$. Expansion of the displayed
polynomials gives
\[
 \begin{aligned}
 c_0&=\frac{370815526253}{62986315392000}>\frac1{200},\\
 c_1&=-\frac{697026734017}{22230464256000}>-\frac1{30},\qquad
 c_2=\frac{426624588209}{1587890304000}>\frac4{15}.
 \end{aligned}
\]
Among $c_3,\ldots,c_{17}$, precisely $c_3,c_5,c_6$ are positive.
The same expansion therefore gives
\[
 \sum_{j=3}^{17}|c_j|
 =2(c_3+c_5+c_6)-\bigl(\mathcal R(2)-c_0-c_1-c_2\bigr)
 =\frac{35569199059327}{377917892352000}<\frac1{10}.
\]
Hence for $|y|\le1$,
\[
 \begin{aligned}
 \mathcal R(1+y)
 &>\frac1{200}-\frac{|y|}{30}
    +\left(\frac4{15}-\frac1{10}\right)y^2\\
 &=\frac1{300}+\frac16\left(|y|-\frac1{10}\right)^2>0.
 \end{aligned}
\]
Thus $G(s,2)>57/100$ for $0\le s\le2$.

Set $s=m/M$, $u=N/M$, and $S_j=M^{-1}\sum_{i=1}^j f_2(i/M)$.
Differentiating the five factorial prefixes gives exactly
\[
 V_A/M^3=S_m+S_N+S_{m+N}-S_{2m}-S_{2N}.
\]
Since $f_2$ decreases and $f_2(0)=1$, comparison with the left and
right endpoint sums on the same partition gives
\[
 0\le F_2(j/M)-S_j\le\frac{f_2(0)-f_2(j/M)}M\le\frac1M.
\]
Use the lower bounds for the three positive prefixes and the upper
bounds for the two prefixes with a minus sign. This gives
$V_A/M^3\ge G(s,u)-3/M$. Finally,
$V_A-V_B=M^2 f_2(1/M)\le M^2$, so $V_B/M^3>0.57-4/M$.
\end{proof}

\subsection{The positive recurrence on the outer arc}
We realize $\qbinom{L}{K}_{qe^{i\theta}}/\qbinom{L}{K}_q$, up to a phase,
as the conditional characteristic function of a sum of selected
indices. We first estimate that characteristic function for
independent Bernoulli trials conditioned on their total.

\begin{lemma}\label{lem:conditioning}
Let $X_1,\ldots,X_L$ be independent Bernoulli variables with parameters
$0<p_j<1$. Suppose $K$ is a mode of $T=\sum_jX_j$. Put
\[
 w_j=p_j(1-p_j),\qquad W=\sum_jw_j,\qquad
 R=\left|\sum_jw_je^{ij\theta}\right|.
\]
Then
\[
 \left|\mathbb E\left[e^{i\theta\sum_jjX_j}\mid T=K\right]\right|
 \le24e^{-19(W-R)/20}.
\]
\end{lemma}
\begin{proof}
Chebyshev's inequality places at least $3/4$ of the mass in an interval
containing at most $4\sqrt W+1$ integers. Since $K$ is a mode,
\[
 \mathbb P(T=K)\ge\frac{3}{4(4\sqrt W+1)}
 \ge\frac1{6(1+\sqrt W)}.
\]
Write $Y=\sum_j jX_j$ and $\beta=\arg\sum_jw_je^{ij\theta}$,
with $\beta=0$ when $R=0$. Fourier inversion gives
\[
 \mathbb E[e^{i\theta Y}\mathbf{1}_{\{T=K\}}]
 =\frac1{2\pi}\int_{-\pi}^{\pi}e^{-iKt}
      \prod_{j=1}^L(1-p_j+p_je^{i(t+j\theta)})\,dt .
\]
The square of the modulus of a single factor is
$1-2w_j(1-\cos(t+j\theta))$. Using $1-x\le e^{-x}$ and rotating
the angle by $\beta$, we therefore bound the conditional modulus by
\[
 6(1+\sqrt W)e^{-(W-R)}
 \frac1{2\pi}\int_{-\pi}^{\pi} e^{-R(1-\cos t)}\,dt .
\]
The integral divided by $2\pi$ is at most one. For $R>0$, the bound
$1-\cos t\ge2t^2/\pi^2$ gives the additional upper bound
$\sqrt\pi/(2\sqrt{2R})$ by integration over the real line.
If $R\ge W/2$, the prefactor is at most $12$ when $W<1$, and at
most $6\sqrt\pi<11$ when $W\ge1$. If $R<W/2$, absorb one twentieth
of the exponential cost and use
\[
 6(1+\sqrt W)e^{-(W-R)/20}
 \le6(1+\sqrt W)e^{-W/40}
 \le6(1+\sqrt{20/e})<24.
\]
\end{proof}

For a Gaussian polynomial with $1\le K<L$, take
\[
 z=q^{-K}\frac{[K]_q}{[L-K+1]_q},\qquad
 p_j=\frac{zq^j}{1+zq^j}.
\]
By expanding $\prod_{j=1}^L(1+zq^j t)$, the probability of $T=r$
is proportional to $z^rq^{r(r+1)/2}\qbinom{L}{r}_q$. Thus
$\mathbb P(T=r)/\mathbb P(T=r-1)=zq^r[L-r+1]_q/[r]_q$.
For $0<q\le1$ this decreases in $r$ and is one at $r=K$, so
$K-1$ and $K$ are modes. Conditional on $T=K$, the weight of a set
$S$ is proportional to $q^{\sum_{j\in S}j}$. Applying the same
product identity at $qe^{i\theta}$ and at $q$ gives
\[
 \mathbb E\!\left[e^{i\theta\sum_{j=1}^L jX_j}\mid T=K\right]
 = e^{iK(K+1)\theta/2}\,\frac{\qbinom{L}{K}_{qe^{i\theta}}}{\qbinom{L}{K}_q}.
\]

For larger angles we use the positive recurrence
\eqref{eq:wz-mixture}. Its Gaussian polynomial and its inner kernel
have different dependence on the summation index $h$. We first estimate
the Gaussian polynomial, using the conditional Bernoulli representation
of Lemma~\ref{lem:conditioning}.
\begin{lemma}\label{lem:gaussian-entropy}
Let $N\ge500$ and let $K,L$ be integers with
$1\le L\le N$, $0\le K\le L$, and
$0\le\lambda N\le2$. Put $q=e^{-\lambda}$,
$v=K(L-K)/L$, and $\rho_\theta=\csc(|\theta|/2)$.
Then, for $0<|\theta|\le\pi$,
\[
 \left|\frac{\qbinom{L}{K}_{qe^{i\theta}}}{\qbinom{L}{K}_q}\right|
 \le31\exp\left\{-\frac{171}{200}v+
                         \frac{19}{80}\rho_\theta\right\}.
\]
\end{lemma}
\begin{proof}
When $K=0$ or $K=L$, the Gaussian polynomial equals one and $v=0$.
For $1\le K<L$, use the fugacity that makes $K$ a mode,
\[
 z=q^{-K}\frac{[K]_q}{[L-K+1]_q},\qquad
 p(x)=\frac{ze^{-\lambda x}}{1+ze^{-\lambda x}},\quad 0\le x\le L+1.
\]
At $\lambda=0$, use the continuous limits
$z=K/(L-K+1)$ and $p(x)=K/(L+1)$.
The exact identity
\[
 \frac{1+z}{1+ze^{-\lambda(L+1)}}=e^{\lambda K}
\]
shows that $\int_0^{L+1}p(x)\,dx=K$.
Set $T=L+1$, $S=\lambda T$, $x=\lambda K$, and $y=\lambda(T-K)$.
For $\lambda>0$, the identity $p'=-\lambda p(1-p)$ gives
\[
 W_c:=\int_0^T p(t)(1-p(t))\,dt
 =\frac{p(0)-p(T)}{\lambda}
 =\frac{2\sinh(x/2)\sinh(y/2)}{\lambda\sinh(S/2)}.
\]
Consequently
\[
 \frac{W_c}{K(1-K/T)}
 =\frac{\operatorname{sinhc}(x/2)\operatorname{sinhc}(y/2)}
        {\operatorname{sinhc}(S/2)}
 \ge\frac{\operatorname{sinhc}(S/4)^2}{\operatorname{sinhc}(S/2)}
 =\frac{\tanh(S/4)}{S/4}.
\]
Here $\operatorname{sinhc}t=\sinh(t)/t$, continuously extended at zero.
The inequality follows because
$(\log\operatorname{sinhc}t)''=t^{-2}-\sinh^{-2}t\ge0$.
Since $S\le2(1+1/500)=501/250$ and
$\tanh t\ge t-t^3/3$, the final ratio is greater than $9/10$.
The same conclusion holds at $\lambda=0$ by continuity.

The function $w(t)=p(t)(1-p(t))$ is unimodal and bounded by $1/4$.
Each of its superlevel sets in $[0,L+1]$ is an interval. The number of
integers from $1$ to $L$ in an interval of length $\ell$ is at least
$\ell-1$. Integrating this inequality over the levels of $w$ gives
\[
 W:=\sum_{j=1}^Lw(j)\ge W_c-\frac14
 \ge\frac9{10}K(1-K/(L+1))-\frac14
 \ge\frac9{10}v-\frac14.
\]
Abel summation and unimodality give
\[
 R(\theta):=\left|\sum_{j=1}^Lw(j)e^{ij\theta}\right|
 \le\frac{w(1)+w(L)+\sum_{j=1}^{L-1}|w(j+1)-w(j)|}
              {2|\sin(\theta/2)|}
 =\frac{\max_jw(j)}{|\sin(\theta/2)|}
 \le\frac14\rho_\theta.
\]
Lemma~\ref{lem:conditioning} now bounds
the normalized modulus by $24\exp[-19(W-R)/20]$, hence by
\[
 24e^{19/80}
 \exp\left\{-\frac{171}{200}v+\frac{19}{80}\rho_\theta\right\}.
\]
To absorb the constant $e^{19/80}$, use the exponential series with
$t=19/80$. After its quadratic term the successive ratios are at most
$t/3$, giving
\[
 e^t\le1+t+\frac{t^2}{2(1-t/3)}
 =\frac{44841}{35360}<\frac{31}{24}.
\]
This proves the lemma.
\end{proof}

Alongside the Gaussian polynomial, \eqref{eq:wz-mixture} contains
$H(b,r;q)$ with $r=p-2h$. A choice of Bernoulli parameters centers its quadratic
exponent at their mean. Parseval's identity then gives decay in $r$.
\begin{lemma}\label{lem:inner-parseval}
Suppose $N\ge500$, $0\le\lambda N\le2$, and $q=e^{-\lambda}$.
Let $r,b$ be integers with $0\le r\le b\le N$. For $0<|\theta|\le\pi$,
\[
 \frac{|H(b,r;qe^{i\theta})|}{H(b,r;q)}
 \le\frac53\sqrt{r+1}
 \exp\left\{-\frac{13r}{125}+\frac14\csc(|\theta|/2)\right\}.
\]
\end{lemma}
\begin{proof}
For $r=0$, $H(b,0;q)=1$ and the asserted upper bound exceeds one.
For $r>0$, there is a unique
$\nu\in(0,r/2]$ satisfying
\[
 \nu=\sum_{j=1}^r\frac1{1+\exp\{\lambda(b-r+j-1/2+\nu)\}}.
\]
The right-hand side is continuous and nonincreasing, is positive at
zero, and is at most $r/2$ at $r/2$. The left-hand side increases
strictly. This proves existence and uniqueness, including $\lambda=0$.
Put
\[
 z=q^{b-r-1/2+\nu},\qquad
 p_j=\frac{zq^j}{1+zq^j},\quad
 Z=\prod_{j=1}^r(1+zq^j),\quad
 W=\sum_jp_j(1-p_j).
\]
Let $J$ be the sum of independent Bernoulli variables with success
probabilities $p_j$. Thus $\mathbb EJ=\nu$ and $\operatorname{Var}J=W$.
The coefficient of $w^j$ in $\prod_{\ell=1}^r(1+zq^\ell w)$ is
$z^j q^{j(j+1)/2}\qbinom{r}{j}_q$.
The choice of $z$ gives
\[
 \frac{q^{j(b-r+j)}}{z^j q^{j(j+1)/2}}
 =q^{-\nu^2/2}q^{(j-\nu)^2/2}.
\]
Taking expectation and using Jensen's inequality yields
\[
 H(b,r;q)=q^{-\nu^2/2}Z\,
       \mathbb E\exp[-\lambda(J-\nu)^2/2]
 \ge q^{-\nu^2/2}Z e^{-\lambda W/2}.
\]
If $C_j$ are the coefficients of
$\prod_{\ell=1}^r(1+zq^\ell e^{i\ell\theta}w)$, then
\[
 |H(b,r;qe^{i\theta})|
 \le q^{-\nu^2/2}\sum_{j=0}^r
       q^{(j-\nu)^2/2}|C_j|
 \le q^{-\nu^2/2}\sqrt{r+1}\left(\sum_j|C_j|^2\right)^{1/2}.
\]
Set $w_j=p_j(1-p_j)$ and $R=|\sum_jw_je^{ij\theta}|$. Parseval
gives
\[
 \frac{\sum_j|C_j|^2}{Z^2}
 =\frac1{2\pi}\int_{-\pi}^\pi
   \prod_{j=1}^r |1-p_j+p_je^{i(t+j\theta)}|^2\,dt
 \le e^{-2W+2R},
\]
because $|1-p+pe^{iv}|^2\le e^{-2p(1-p)(1-\cos v)}$.
Combining the two kernel estimates gives
\[
 \frac{|H(b,r;qe^{i\theta})|}{H(b,r;q)}
 \le e^{\lambda W/2}\sqrt{r+1}\,e^{-W+R}.
\]
The sequence $w_j$ is unimodal, so Abel summation again gives
$R\le\tfrac14\csc(|\theta|/2)$. Also $W\le r/4$ and $\lambda N\le2$,
so $\lambda W/2\le1/4$.

We claim $W\ge(13/125)r-1/4$. Define
$w(t)=1/(4\cosh^2(t/2))$. Since $0<p_j\le1/2$ and
$\nu\le r/2$, setting $h=r/N\le1$ gives
\[
 p_j(1-p_j)\ge w\left(2-h+\frac{2j-1}{N}\right).
\]
These arguments are the midpoints of $r$ intervals of length $2/N$.
The level-set counting bound used for the Gaussian polynomial gives
\[
 W\ge\frac N2\int_{2-h}^{2+h}w(t)\,dt-\frac14.
\]
The symmetric average of $w$ on this interval is
\[
 \frac1{2h}\int_{2-h}^{2+h}w(t)\,dt
 =\frac{\sinh h}{2h(\cosh2+\cosh h)}
 \ge\frac1{2(\cosh2+1)}>\frac{13}{125}.
\]
The middle inequality is equivalent to
$\cosh2(\sinh h-h)\ge h\cosh h-\sinh h$.
Since $\cosh2>3$, it suffices to show $4\sinh h-h\cosh h-3h\ge0$.
This function and its first derivative vanish at zero. Its second
derivative is $2\sinh h-h\cosh h\ge0$ on $[0,1]$, since
$\tanh h\ge h-h^3/3\ge h/2$. The final numerical inequality follows from
$\cosh2<99/26$, which follows, for example, from $e^2<37/5$ and
$e^{-2}<1/7$. Thus the claim is proved.
Combining it with the preceding estimate gives the result, since
$e^{1/2}<5/3$.
\end{proof}

We next compare the weights in \eqref{eq:wz-mixture} at $q$ and at
$1$. This permits an estimate of the cost of large $h$ using the
untilted occupancy distribution.
\begin{lemma}\label{lem:weight-domination}
For fixed integers $a\ge0$, $p\ge0$, $b=a+p$, let $W_h(q)$ be the
normalized weights of \eqref{eq:wz-mixture} at $0<q\le1$.
The distribution $W(q)$ is stochastically dominated by $W(1)$.
Equivalently, its expectation of every nondecreasing function of $h$
is at most the corresponding expectation at $q=1$.
\end{lemma}
\begin{proof}
For $h<\lfloor p/2\rfloor$ put $r=p-2h\ge2$. Applying
\eqref{eq:H-ratio} twice in \eqref{eq:tilted-weight-ratio} gives
\[
 \frac{W_{h+1}(q)}{W_h(q)}
 \le\frac{q^{2a+4h+2}}{(1+q^b)^2}
     \frac{[r]_q[r-1]_q}{[h+1]_q[a+h+1]_q}.
\]
Put $q=e^{-\lambda}$, using $\operatorname{sinhc}(0)=1$ at $q=1$.
The right-hand side equals
\[
 \frac{r(r-1)}{4(h+1)(a+h+1)}
 \frac{e^{-\lambda(a/2+h+1/2)}}{\cosh^2(\lambda b/2)}
 \frac{\operatorname{sinhc}(\lambda r/2)
       \operatorname{sinhc}(\lambda(r-1)/2)}
      {\operatorname{sinhc}(\lambda(h+1)/2)
       \operatorname{sinhc}(\lambda(a+h+1)/2)}.
\]
Since $\operatorname{sinhc}x\ge1$ and
$\operatorname{sinhc}x\le\cosh x$ for $x\ge0$, while $r\le b$,
the final two factors have product at most one. We therefore obtain
\[
 \frac{W_{h+1}(q)}{W_h(q)}
 \le\frac{r(r-1)}{4(h+1)(a+h+1)}
 =\frac{W_{h+1}(1)}{W_h(1)}.
\]
Thus $W_h(q)/W_h(1)$ is nonincreasing and its $W(1)$-mean is one.
The difference $W_h(q)-W_h(1)$ changes sign at most once, from positive
to negative. Every upper tail under $W(q)$ is therefore bounded by its
counterpart under $W(1)$. This includes one-point supports.
\end{proof}

The bounds for the Gaussian factor and for $H$ contribute a quadratic
function $\mathcal J_p(h)$ to the exponent. We need its exponential average
under the weights in \eqref{eq:wz-mixture}.

\begin{lemma}\label{lem:laplace-cost}
Let $p\in\Z$, $p\ge100$, $\gamma=171/200$, and $0\le\delta\le1/4$. For
\[
 \mathcal J_p(h)=\delta p+2(\gamma-\delta)h-4\gamma h^2/p,
\]
one has, uniformly over the weights in \eqref{eq:wz-mixture} at $0<q\le1$,
\[
 \mathbb E_q e^{-\mathcal J_p(h)}\le2e^{-\delta p}.
\]
\end{lemma}
\begin{proof}
For $j\ge0$ write $(n)_j=n(n-1)\cdots(n-j+1)$, with $(n)_0=1$.
To count $j$ distinguished fully selected pairs in the uniform
$p$-subset, first choose them from the $b$ pairs and then require all
their $2j$ members to be selected. This gives, for $0\le2j\le p$,
\[
 \mathbb E_1\binom{h}{j}
 =\binom bj\frac{(p)_{2j}}{(2b)_{2j}}
 =\frac{(p)_{2j}}{2^j j!\prod_{i=0}^{j-1}(2b-2i-1)}.
\]
These quantities decrease with $b\ge p$. For $t\ge1$, expand
$t^h=\sum_j\binom{h}{j}(t-1)^j$. The coefficients are nonnegative, so the
untilted moment $\mathbb E_1 t^h$ is largest at $b=p$.
When $b=p$, the generating function by subset size and fully selected
pairs gives
\[
 \mathbb E_1 t^h=\frac{[x^p](1+2x+tx^2)^p}{\binom{2p}{p}}
 \le\frac{(2+2\sqrt t)^p}{\binom{2p}{p}}
 \le2\sqrt p\left(\frac{1+\sqrt t}{2}\right)^p .
\]
For the first inequality evaluate the polynomial at $x=t^{-1/2}$
and divide by $x^p$. For the second,
$c_p=\sqrt p\binom{2p}{p}/4^p$ satisfies $c_1=1/2$ and
$(c_{p+1}/c_p)^2=1+1/(4p(p+1))>1$. Stochastic domination transfers
the moment bound to $\mathbb E_q t^h$. For $1/2\le z<1$, Markov's
inequality with $\sqrt t=z/(1-z)$ now gives
\[
 \mathbb P_q(h\ge pz/2)\le 2\sqrt p\,e^{-pI(z)},
 \qquad I(z)=\log 2+z\log z+(1-z)\log(1-z).
\]
At $z=1$, take $t\to\infty$ and set $0\log0=0$.

Put $z=2h/p$. For $z\le7/10$,
\[
 \mathcal J_p(h)/p=\delta+(\gamma-\delta)z-\gamma z^2\ge\delta,
\]
since $3\gamma/10>1/4\ge\delta$.
For $z\ge7/10$, define
$\mathcal K(z)=I(z)+\delta+(\gamma-\delta)z-\gamma z^2$.
Here $\mathcal K''=1/(z(1-z))-2\gamma>0$ and
\[
 \mathcal K'(7/10)\ge4(7/10-1/2)-2\gamma/5-\delta>0.
\]
Moreover $I(7/10)>2(7/10-1/2)^2=2/25$, because $I''\ge4$.
It follows that
\[
 \mathcal K(z)\ge \mathcal K(7/10)>\delta+\frac2{25},
\]
using $21\gamma/100-7\delta/10\ge91/20000>0$.
The small-$z$ terms contribute at most $e^{-\delta p}$.
For each remaining $h$, bound its point probability by the upper-tail
probability at $z=2h/p$. There are at most $p/2+1$ such integers, so
their total contribution is at most
\[
 (p+2)\sqrt p\,e^{-(\delta+2/25)p}\le e^{-\delta p}.
\]
Indeed, the logarithmic derivative of $(p+2)\sqrt p\,e^{-2p/25}$ is
$1/(p+2)+1/(2p)-2/25<0$ for $p\ge100$, and at $100$ the value is
$1020e^{-8}<1$.
This proves the lemma.
\end{proof}

Averaging the two factor estimates yields the required outer-arc bound.
\begin{lemma}\label{lem:wz-outer-bound}
Let $a,b\in\Z_{\ge0}$, $N=\max(a,b)\ge500$, $q=e^{-\lambda}$, and
$0\le\lambda N\le2$. For $0<|\theta|\le\pi$,
\[
 \left|\frac{A_{a,b}(qe^{i\theta})}{A_{a,b}(q)}\right|
 \le\min\left\{1,C\sqrt{N+1}
  \exp\left(-cN+\frac{39}{80}\csc(|\theta|/2)\right)\right\},
\]
with the constants
\[
 \begin{array}{c|c|c}
 \lambda N&C&c\\\hline
 {[0,2]}&104&26/625\\
 {[0,1/4]}&64&12/125.
 \end{array}
\]
\end{lemma}
\begin{proof}
Keep $N$ and $\lambda$ fixed in each application of the recurrence.
Write
$\rho_\theta=\csc(|\theta|/2)$ and $\gamma=171/200$. Suppose the inner
kernel satisfies, uniformly for $0\le r\le b\le N$,
\[
 \frac{|H(b,r;qe^{i\theta})|}{H(b,r;q)}
 \le C_H\sqrt{r+1}\,e^{-\delta r+\rho_\theta/4},
 \qquad 0\le\delta\le1/4 .
\]
For a subproblem with $a\le b\le N$ and $p=b-a\ge100$, normalize
\eqref{eq:wz-mixture} at $q$. The uniform bound one for the child,
Lemma~\ref{lem:gaussian-entropy} with $(L,K)=(p,2h)$,
and the preceding inner bound for $r=p-2h$ give
\[
 \frac{|A_{a,b}(qe^{i\theta})|}{A_{a,b}(q)}
 \le31C_H\sqrt{N+1}\,e^{39\rho_\theta/80}\,
       \mathbb E_q e^{-\mathcal J_p(h)},
 \qquad \mathcal J_p(h)=\delta(p-2h)+2\gamma h(1-2h/p).
\]
Lemma~\ref{lem:laplace-cost} bounds the expectation by $2e^{-\delta p}$.

For the original polynomial take $b=N$ by symmetry. If $p\ge2N/5$,
this gives an exponential cost of at least $2\delta N/5$. Otherwise,
in every summand of the first recurrence the child $A_{a,h}=A_{h,a}$
has gap $a-h=N-p-h\ge N-3p/2>2N/5$. Apply the same estimate to each
child, bound the remaining normalized factors by one, and average with
the positive weights. Every gap to which the estimate is applied is
at least $2N/5\ge200$. In both cases,
\begin{equation}\label{eq:common-wz-far-bound}
 \frac{|A_{a,b}(qe^{i\theta})|}{A_{a,b}(q)}
 \le\min\{1,62 C_H\sqrt{N+1}\,
                 e^{-2\delta N/5+39\rho_\theta/80}\}.
\end{equation}
For $\lambda N\le2$, Lemma~\ref{lem:inner-parseval}
supplies $C_H=5/3$, $\delta=13/125$. Since $62(5/3)<104$ and
$2(13/125)/5=26/625$, the claimed bound follows.

For the second row, suppose $\lambda N\le1/4$.

For $r=0$ the inner kernel is one and its bound below holds directly.
For $r\ge1$, in the Bernoulli representation of $H(b,r;q)$
the exponent in the denominator of $p_j$ is
\[
 x_j=\log\frac{1-p_j}{p_j}
    =\lambda(b-r+j-1/2+\nu).
\]
Since $1\le j\le r$ and $0<\nu\le r/2$, it satisfies
$0\le x_j\le\lambda(b+r/2)\le3/8$. Hence
\[
 p_j(1-p_j)\ge\frac1{4\cosh^2(3/16)}
 \ge\frac14-\frac{9}{1024}>\frac6{25},
\]
where $\operatorname{sech}^2x=1-\tanh^2x\ge1-x^2$.
Thus $W\ge6r/25$. The same weighted Parseval argument gives
\[
 \frac{|H(b,r;qe^{i\theta})|}{H(b,r;q)}
 \le\frac{32}{31}\sqrt{r+1}
       \exp[-6r/25+\rho_\theta/4],
\]
because $e^{\lambda W/2}\le e^{1/32}\le32/31$.
Apply \eqref{eq:common-wz-far-bound} with $C_H=32/31$ and
$\delta=6/25$. The identities $62(32/31)=64$ and
$2(6/25)/5=12/125$ give the stated constants.
\end{proof}

\subsection{The compact arc and the paired saddle}

Between the first lobe and the outer arc, the five factorials in
$A_{m,N}$ have a continuum limit with an explicit discretization error.
Write $x=m/N$, $q=e^{-u/N}$, $t=N\theta$, and
\[
 A_{m,N}(q)=\frac{(q;q)_{2m}(q;q)_{2N}}
 {(q;q)_m(q;q)_N(q;q)_{m+N}},\qquad
 \phi_A(\theta)=\frac{A_{m,N}(qe^{i\theta})}{A_{m,N}(q)}.
\]

\begin{lemma}\label{lem:compact-arc}
For $0\le m\le N$ and $1/4\le u\le2$,
\begin{equation}\label{eq:compact-arc}
 N\ge1000,\quad 3\le N\theta\le100
 \quad\Longrightarrow\quad |\phi_A(\theta)|\le e^{-N/16}.
\end{equation}
\end{lemma}
\begin{proof}
We first obtain a uniform bound for the polynomial in the Lambert
expansion of $\log A_{m,N}$. Let
\[
 Q(a,b)=1-a-b-ab+a^2+b^2=(1-a)(1-b)+(a-b)^2.
\]
Then $|Q(a,b)|\le4$ on the closed bidisc. By the maximum modulus theorem,
it suffices to put $a=e^{i(s+v)}$, $b=e^{i(s-v)}$. Put $c=\cos v$ and
$d=4c^2-3$; replacing $s$ by $s+\pi$ absorbs a change of sign of $c$. Then
$Q=1-2ce^{is}+de^{2is}$, and, with $X=\cos s$,
\[
 |Q|^2=(1-d)^2+4c^2-4c(1+d)X+4dX^2.
\]
At $X=\pm1$, its square root is at most $|4c^2-2|+2c\le4$.
If $d\ge0$, convexity proves the assertion. If $d<0$ and the quadratic
has its maximum inside $[-1,1]$, its vertex is
$X=c(1+d)/(2d)$. In this case $d\le-1/4$, since for $-1/4<d<0$ one has
$c>3/4$ and $1+d>3/4$, so $|X|>1$. At the vertex,
\[
 16-|Q|^2=-\frac{(d+3)(3d^2-18d-1)}{4d}\ge0
 \qquad(-3\le d\le-1/4).
\]
This proves the bound.

Let $\operatorname{Li}_2(z)=\sum_{\ell\ge1}z^\ell/\ell^2$ for $|z|\le1$.
For $\Re w>0$, put
\begin{align*}
 \mathcal C_x(w)&=\frac{1}{w}\sum_{\ell\ge1}
 \frac{Q(e^{-xw\ell},e^{-w\ell})}{\ell^2}\\
 &=\frac{1}{w}\biggl\{\frac{\pi^2}{6}
 -\operatorname{Li}_2(e^{-xw})-\operatorname{Li}_2(e^{-w})
 -\operatorname{Li}_2(e^{-(1+x)w})
 +\operatorname{Li}_2(e^{-2xw})+\operatorname{Li}_2(e^{-2w})\biggr\}.
\end{align*}
The first expression and the polynomial bound give
\begin{equation}\label{eq:continuum-modulus}
 |\mathcal C_x(w)|\le\frac{2\pi^2}{3|w|}<\frac{20}{3|w|}.
\end{equation}
Take the analytic logarithm of the product for $|z|<1$ that vanishes
at $z=0$. For integer sequences $N\to\infty$, $m/N\to x$, the exact Lambert expansion
\[
 \log A_{m,N}(z)=\sum_{\ell\ge1}
 \frac{z^\ell Q(z^{m\ell},z^{N\ell})}{\ell(1-z^\ell)}
\]
and dominated convergence show that
$N^{-1}\log A_{m,N}(e^{-w/N})\to \mathcal C_x(w)$. Domination follows from
$|z^\ell/[N(1-z^\ell)]|\le1/(\ell\Re w)$ and $|Q|\le4$.
In particular, known coefficient nonnegativity implies
$\Re \mathcal C_x(u-it)\le \mathcal C_x(u)$ for $u>0$. On the real axis the canceled product is bounded by the partition
product. Expanding the logarithm of that product and using
$e^{u\ell/N}-1\ge u\ell/N$ gives
\[
 \frac1N\log A_{m,N}(e^{-u/N})
 \le \frac1N\sum_{\ell\ge1}\frac{1}{\ell(e^{u\ell/N}-1)}
 \le \frac{\pi^2}{6u}.
\]
The same bound holds for $\mathcal C_x(u)$.

We next bound $\mathcal C_x(u)$ from below, uniformly in $x$. Each surviving factor has the form $(1-q^a)/(1-q^b)$ with $a>b$.
Its logarithm decreases with $u$ when $q=e^{-u/N}$, because
$a/(e^{ua/N}-1)<b/(e^{ub/N}-1)$. Thus $\mathcal C_x(u)\ge \mathcal C_x(2)$ for $u\le2$.
For $0<x<1$, write $y=e^{-2x}$ and $c=e^{-2}$. Differentiation gives
\[
 \frac{d}{dx}\mathcal C_x(2)=\log\frac{(1-y)(1+y)^2}{1-cy}.
\]
Its sign is that of $(1-y)(1+y)^2-(1-cy)=y(1+c-y-y^2)$.
As $x$ increases, $y$ decreases from $1$ to $c$, and
$1+c-y-y^2$ increases from $c-1<0$ to $1-c^2>0$.
Continuity at both endpoints now gives the unique minimum at
\[
 y_*=(\sqrt{5+4e^{-2}}-1)/2,\qquad x_*=-\tfrac12\log y_*.
\]
Put
\[
 c_-:=\frac{13533}{100000},\qquad c_+:=\frac{13534}{100000},\qquad
 y_-:=\frac{677}{1000},\qquad y_+:=\frac{67701}{100000}.
\]
The exponential series gives
\[
 \frac1{c_+}<\sum_{j=0}^{10}\frac{2^j}{j!}<e^2
 <\sum_{j=0}^{10}\frac{2^j}{j!}
       +\frac65\frac{2^{11}}{11!}<\frac1{c_-}.
\]
For the upper bound, the ratio of consecutive terms in the tail
starting with $2^{11}/11!$ is at most $1/6$.
The rational comparisons
\[
 y_-^2+y_-<1+c_-,\qquad y_+^2+y_+>1+c_+
\]
therefore imply
\[
 c_-<c< c_+,\qquad y_-<y_*<y_+.
\]
For $0<z<1$, positivity of the dilogarithm series gives
\[
 \sum_{j=1}^{10}\frac{z^j}{j^2}
 <\operatorname{Li}_2(z)
 <\sum_{j=1}^{10}\frac{z^j}{j^2}
       +\frac{z^{11}}{121(1-z)}.
\]
At the four displayed rational endpoints these finite sums satisfy
\[
 \sum_{j=1}^{10}\frac{y_-^{2j}+c_-^{2j}}{j^2}
 >\frac{54388}{100000},
\]
\[
 \sum_{z\in\{y_+,c_+,y_+c_+\}}
 \left(\sum_{j=1}^{10}\frac{z^j}{j^2}
       +\frac{z^{11}}{121(1-z)}\right)
 <\frac{108455}{100000}.
\]
The elementary bound $\pi>314159/100000$ yields
$\pi^2/6>164493/100000$.
Substitution in the formula for $\mathcal C_{x_*}(2)$ now gives
\[
 \begin{aligned}
 2\mathcal C_{x_*}(2)
 &=\frac{\pi^2}{6}
   +\operatorname{Li}_2(y_*^2)+\operatorname{Li}_2(c^2)
   -\operatorname{Li}_2(y_*)-\operatorname{Li}_2(c)
   -\operatorname{Li}_2(cy_*)\\
 &>\frac{164493+54388-108455}{100000}
  =\frac{110426}{100000}.
 \end{aligned}
\]
Consequently
\begin{equation}\label{eq:real-entropy}
 \mathcal C_x(u)>\frac{55213}{100000}>\frac{69}{125}
 \qquad(0\le x\le1,\quad0<u\le2).
\end{equation}

\smallskip\noindent\emph{The harmonic boundary estimate.}
For the complex tilt we claim
\begin{equation}\label{eq:harmonic-boundary}
 \Re \mathcal C_x(u-it)<\frac{12}{25}=0.48
 \qquad(0\le x\le1,\ u\ge0,\ t\ge3).
\end{equation}
For fixed $x$, the real part is harmonic in $u>0,t>3$.
The series defining $w\mathcal C_x(w)$ converges uniformly for $\Re w\ge0$,
since $|Q|\le4$ and $\sum_{\ell\ge1}\ell^{-2}<\infty$.
It thus extends continuously to the two boundary rays, where $|w|\ge3$.
By \eqref{eq:continuum-modulus} the function tends uniformly to zero at infinity.
Apply the maximum principle first on the intersection with $|w|<R$
and then let $R\to\infty$. It suffices to bound the two rays.
Write $\operatorname{Cl}_2(v)=\sum_{\ell\ge1}\sin(\ell v)/\ell^2$.
Taking the real part on $u=0$ gives
\begin{align*}
 t\Re \mathcal C_x(-it)={}&\operatorname{Cl}_2(xt)+\operatorname{Cl}_2(t)
 +\operatorname{Cl}_2((1+x)t)\\
 &-\operatorname{Cl}_2(2xt)-\operatorname{Cl}_2(2t).
\end{align*}
For $t\ge18$, the elementary bound
$|\operatorname{Cl}_2(v)|\le\pi^2/6<5/3$ gives an upper bound $25/54<12/25$.
On $t=3,u\ge4$, coefficient nonnegativity and the partition bound give
$\Re \mathcal C_x(u-3i)\le \mathcal C_x(u)\le\pi^2/(6u)<5/12<12/25$.
The remaining two rectangles are
\[
 (x,t)\in[0,1]\times[3,18],\ u=0;
 \qquad (x,u)\in[0,1]\times[0,4],\ t=3.
\]
Set $L=32$ and truncate each dilogarithm at order $L$.
The five tails together have absolute value at most $5/(L|w|)$.
For the two rectangles above, Python integer interval arithmetic
encloses the truncated sums and their derivatives. Exponential and
trigonometric evaluations use Taylor remainders bounded by
$3/43!$ and $4^{62}/62!$ after range reduction; both are less than
$10^{-36}$. Machin's identity, evaluated by the alternating
arctangent series in the same interval arithmetic, encloses $\pi$.

At each cell centre we enclose the truncated real part $F$
and both first derivatives. With $H_L=\sum_{\ell\le L}1/\ell$, bounds
for the absolute values of its second derivatives are
\[
\begin{array}{c|ccc}
 & |F_{xx}|&|F_{xy}|&|F_{yy}|\\ \hline
 u=0,\ y=t &6tL&7L&14L/3+16H_L/9+50/81\\
 t=3,\ y=u &30L&7L&14L/3+16H_L/9+20/27.
\end{array}
\]
To verify the table, write $a=(x,1,1+x,2x,2)$ and $b=da/dx$.
Then
\[
 \sum_j b_j^2=6,\quad \sum_j a_jb_j\le7,\quad
 \sum_j a_j^2\le14,\quad \sum_j a_j\le8.
\]
Differentiating the five exponential sums gives
$|F_{xx}|\le |w|L\sum_j b_j^2$ and
$|F_{xy}|\le L\sum_j a_jb_j$. For the $yy$ derivative of the
quotient, the three terms have absolute values at most
$14L/|w|$, $16H_L/|w|^2$ and $2S_0/|w|^3$.
On $u=0$ the numerator is a signed sum of five
truncated Clausen sums, so $S_0=5\pi^2/6$ suffices. On $t=3$
it also contains the exact constant $\pi^2/6$, so $S_0=\pi^2$
suffices. Using $|w|\ge3$, $\pi^2<10$, and
$|w|\le5$ for $t=3$, $0\le u\le4$, gives the displayed entries.
Let $H_{xx},H_{xy},H_{yy}$ be the respective bounds on the whole
cell; in particular, the first entry uses the upper $t$ endpoint.
If the cell has half-widths $h_x,h_y$, its certified upper bound is
\[
 F_c+|F_{x,c}|h_x+|F_{y,c}|h_y+\tfrac12 H_{xx}h_x^2
 +H_{xy}h_xh_y+\tfrac12 H_{yy}h_y^2+\frac{5}{Ld_0}.
\]
Here $d_0$ is the lower $t$ endpoint of a Clausen cell, and $d_0=3$
for a vertical cell, so $|w|\ge d_0$ throughout the cell.
Successive bisection gives partitions on which this upper bound
is strictly less than $12/25$. The exact interval calculation gives
\[
 \begin{array}{c|r|c}
 \text{rectangle}&\text{number of cells}&\text{upper bound}\\\hline
 u=0,\ 0\le x\le1,\ 3\le t\le18&12594&0.479996\\
 t=3,\ 0\le x\le1,\ 0\le u\le4&7519&0.4799992.
 \end{array}
\]
Both bounds are below $12/25$. The maximum principle therefore
proves~\eqref{eq:harmonic-boundary}.

\smallskip\noindent\emph{The discretization error.}
To compare the continuum estimate with the finite product, put $w=u-it$ and
\[
 h_w(s)=\log\frac{\sinh(ws/2)}{ws/2},\qquad
 g(s)=\Re h_w(s)-h_u(s),\qquad g(0)=0.
\]
Choose the logarithm continuous from $h_w(0)=0$.
Since $\sinh(ws/2)/(ws/2)=\int_{-1/2}^{1/2}e^{wsr}\,dr$, the
triangle inequality gives $g(s)\le0$. For the lower bound use
$|\sinh((u-it)s/2)|^2=\sinh^2(us/2)+\sin^2(ts/2)$. Together these give
\[
 -\tfrac12\log(1+(t/u)^2)\le g(s)\le0.
\]
For $s>0$ the individual factor satisfies
\[
 \log\frac{|1-e^{-ws}|}{1-e^{-us}}=\log(|w|/u)+g(s).
\]
In the five factorials the constant cancels because
$2m+2N=m+N+(m+N)$. Thus $\log|\phi_A|$ is exactly the signed sum
of the five right-endpoint sums of $g(j/N)$.
The corresponding signed sum of integrals on the intervals of
lengths $2x,2,x,1,1+x$, with signs $+,+,-,-,-$, equals
$\Re \mathcal C_x(w)-\mathcal C_x(u)$. Indeed, one may fix $x=m/N$ and scale
$(m,N)$ by any common positive integer. The normalized sums converge
to these integrals by continuity of $g$, and the Lambert limit gives
the same limit as $\Re \mathcal C_x(w)-\mathcal C_x(u)$.
For $1\le j\le2N$, let $e_j\in\{-1,0,1\}$ be the exponent of
$1-q^j$ in the canceled Euler product. Apply the trapezoidal formula
on each interval $[(j-1)/N,j/N]$ before summing with these signs:
\[
 \sum_{j=1}^{2N}e_jg(j/N)
 =N\sum_{j=1}^{2N}e_j\int_{(j-1)/N}^{j/N}g(s)\,ds
 +\frac12\sum_{j=1}^{2N}e_j\{g(j/N)-g((j-1)/N)\}+R.
\]
The error kernel on each interval has absolute value at most
$1/(8N^2)$. Since $|e_j|\le1$, the total remainder satisfies
\[
 |R|\le\frac1{8N}\sum_{j=1}^{2N}|e_j|
             \int_{(j-1)/N}^{j/N}|g''(s)|\,ds
 \le\frac1{8N}\int_0^2|g''(s)|\,ds.
\]
The first term is $N\{\Re \mathcal C_x(u-it)-\mathcal C_x(u)\}$.
Since $g(0)=0$, telescoping the second term gives
\[
 \frac12\{g(2x)+g(2)-g(x)-g(1)-g(1+x)\}
 \le\frac34\log(1+(t/u)^2).
\]
Consequently
\begin{equation}\label{eq:trapezoid-error}
 \log|\phi_A(\theta)|
 \le N\{\Re \mathcal C_x(u-it)-\mathcal C_x(u)\}
 +\frac34\log(1+(t/u)^2)
 +\frac1{8N}\int_0^2|g''(s)|\,ds.
\end{equation}

We bound the integral in \eqref{eq:trapezoid-error} by
\begin{equation}\label{eq:second-derivative-integral}
 \int_0^2|g''(s)|\,ds
 \le \frac{147\pi}{529}\frac{u^2+t^2}{t}+\frac{t}{\pi}
 +\frac{\pi(u^2+t^2)}{2u}\sum_{k=1}^{K(t)}\frac1{2k-1}
 +\frac{u^2}{6},
\end{equation}
where $K(t)=\lfloor t/\pi+1/2\rfloor$ (including an empty endpoint cell is harmless).
Indeed,
\[
 h_w''(s)=\frac1{s^2}-\frac{w^2}{4\sinh^2(ws/2)}.
\]
To control the removable singularity at zero, use the product for
$\sinh$ in \eqref{eq:sinh-product}. With $d_j=w^2/(4\pi^2j^2)$,
\[
 h_w''(s)=\sum_{j\ge1}2d_j\frac{1-d_js^2}{(1+d_js^2)^2}.
\]
For $0\le s\le\pi/t$ we have $|d_js^2|\le(1+4/9)/4=13/36$.
Hence termwise differentiation is justified by uniform convergence,
and $2\sum_j|d_j|=|w|^2/12$ gives
\[
 |h_w''(s)|\le\frac{|w|^2}{12}
 \frac{1+13/36}{(1-13/36)^2}=\frac{147}{529}|w|^2.
\]
On the cell
$[(2k-1)\pi/t,(2k+1)\pi/t]$, put $s_k=2k\pi/t$ and
$s_0=(2k-1)\pi/t$. Since $|t(s-s_k)/2|\le\pi/2$, the inequality
$|\sin v|\ge2|v|/\pi$ on this interval gives
\[
 |\sinh(ws/2)|^2\ge u^2s_0^2/4+t^2(s-s_k)^2/\pi^2.
\]
Integration over the full real line bounds the corresponding singular term
by $\pi|w|^2/[2u(2k-1)]$. The $1/s^2$ term contributes at most $t/\pi$.
For real $w=u$, $d_j\ge0$ and $|1-d_js^2|/(1+d_js^2)^2\le1$,
so the same product gives $|h_u''(s)|\le u^2/12$ on the full interval.
These observations prove \eqref{eq:second-derivative-integral}.

If $3\le t\le T$, $1/4\le u\le2$, then \eqref{eq:second-derivative-integral} is at most
\[
 J_T=\frac{462}{529}(T+4/3)+\frac T3+
 \frac{11}{7}(1/4+4T^2)\sum_{k=1}^{\lfloor T/3+1/2\rfloor}\frac1{2k-1}
 +\frac23.
\]
This uses $3<\pi<22/7$ and the monotonicity of $u+t^2/u$ for $u\le2<t$.
Substitution of the rational endpoints gives
\[
\begin{array}{c|rrrr}
 T&18&40&70&100\\ \hline
 (3/4)\log(1+16T^2)\text{ is less than}&7&8&9&9\\
 (1/8)J_T\text{ is less than}&482&2853&9827&21466.
\end{array}
\]

For $3\le t\le18$, \eqref{eq:real-entropy}, \eqref{eq:harmonic-boundary}, and
\eqref{eq:trapezoid-error} give
\[
 \log|\phi_A|< -\frac9{125}N+7+\frac{482}{N}\le-\frac N{16}
 \qquad(N\ge1000).
\]
For $18\le t\le40$, $40\le t\le70$, and $70\le t\le100$, use
\eqref{eq:continuum-modulus} with the lower endpoint of the interval and the
corresponding error bounds in the table. The three respective rate lower
bounds are
\[
 \frac{69}{125}-\frac{20}{54},\qquad
 \frac{69}{125}-\frac{20}{120},\qquad
 \frac{69}{125}-\frac{20}{210}.
\]
Each exceeds $1/16+(E+R/1000)/1000$ for its table values $E,R$.
For larger $N$, $E/N+R/N^2$ decreases, proving
\eqref{eq:compact-arc}.
\end{proof}

The preceding estimate requires $u\ge1/4$. For smaller tilts,
reciprocity allows us to transfer a bound on an arc at one positive
tilt to the same arc at zero tilt. The price of the transfer depends
only on an upper bound for the variance.

\begin{lemma}\label{lem:strip-transfer}
Let $A(q)=\sum_{j=0}^Da_jq^j$ be a nonzero polynomial with
nonnegative reciprocal coefficients, and suppose
$V_A(s)\le\sigma^2$ for $s\ge0$. Fix $N>0$,
$0<\varepsilon$, and $0<\theta_0<\pi$. Suppose
\begin{align*}
 \left|\frac{A(e^{-\varepsilon/N+i\theta})}
                   {A(e^{-\varepsilon/N})}\right|&\le \beta
 &&(\theta_0\le\theta\le2\pi-\theta_0),\\
 \left|\frac{A(e^{-s+i\theta_0})}{A(e^{-s})}\right|&\le \gamma
 &&(0\le s\le\varepsilon/N).
\end{align*}
Then, uniformly for $0\le s\le\varepsilon/N$ and
$\theta_0\le\theta\le2\pi-\theta_0$,
\[
 \left|\frac{A(e^{-s+i\theta})}{A(e^{-s})}\right|
 \le\exp\left(\frac{\sigma^2\varepsilon^2}{2N^2}\right)
       \max(\beta,\gamma).
\]
For the super Catalan family with $N\ge500$, the exponential factor
may be replaced by $\exp(\varepsilon^2N/10)$.
\end{lemma}
\begin{proof}
The function $F(z)=e^{Dz/2}A(e^{-z})$ is entire and even. It has real
Taylor coefficients, and $F(z+2\pi i)=(-1)^D F(z)$. On the real axis,
$F(s)\ge F(0)=A(1)$ by pairing reciprocal coefficients. Since
$(\log F)'(0)=0$ and $(\log F)''=V_A$,
\[
 \log\frac{F(s)}{F(0)}
 =\int_0^s(s-t)V_A(t)\,dt\le\frac12\sigma^2s^2.
\]
Consider the rectangle
\[
 -\varepsilon/N\le\Re z\le\varepsilon/N,
 \qquad \theta_0\le\Im z\le2\pi-\theta_0.
\]
On both vertical sides, evenness and conjugation reduce the asserted
bound to the first hypothesis. On both horizontal sides, the same
symmetries and the modulus periodicity reduce it to the second.
Thus every boundary value of $|F(z)|/F(0)$ is at most
$e^{\sigma^2\varepsilon^2/(2N^2)}\max(\beta,\gamma)$. The maximum modulus
principle gives the same bound throughout the rectangle. Finally,
divide $|F(s+i\theta)|$ by $F(s)\ge F(0)$. For the super Catalan specialization,
\eqref{eq:variance-origin} gives
$V_A(s)\le V_A(0)\le(2N^3+N^2)/12<N^3/5$ for $N\ge500$.
\end{proof}

Write $A=A_{m,N}$ and $q=e^{-\lambda}$. For the difference
$a_k-a_{k-1}$ we shall choose the tilt for which $\mu_A(\lambda)=k-1/2$.
The local contribution to the resulting
first-difference integral has scale $\lambda I_G$. Near the midpoint $\lambda$ can be of order $N^{-3}$,
so the phase error must retain its factor of $\lambda$. Reciprocity
supplies this factor through the odd derivatives of the centered
logarithm.

\begin{lemma}\label{lem:paired-local}
Suppose $N\ge500$, $0\le m\le N$, and $0\le\lambda N\le2$. Set
\[
 V=V_A(\lambda),\qquad
 \phi(\theta)=e^{-i\mu_A(\lambda)\theta}\frac{A(qe^{i\theta})}{A(q)},
 \qquad h=\frac74,\qquad I_G=\sqrt{\frac{\pi}{2V}} .
\]
Then $(69/1000)N^3\le V\le N^3/5$. On $|\theta|\le h/N$ write
$\log\phi(\theta)=-V\theta^2/2+R(\theta)+iI(\theta)$, with $R,I$ real
and $\log\phi(0)=0$. One has
\begin{equation}\label{eq:small-local}
 |R(\theta)|\le\frac{13}{4000}N^5\theta^4,
 \qquad |I(\theta)|\le\frac{13}{1000}\lambda N^5|\theta|^3,
 \qquad \Re\log\phi(\theta)\le-\frac{71}{200}V\theta^2.
\end{equation}
In particular,
\begin{align}
 \int_0^{h/N}\Re\phi(\theta)\cos(\theta/2)\,d\theta
 &\ge\left(1-\frac{54}{N}\right)I_G,\label{eq:local-cos}\\
 \left|\int_0^{h/N}\Im\phi(\theta)\sin(\theta/2)\,d\theta\right|
 &\le\frac{10\lambda}{N}I_G.\label{eq:local-sin}
\end{align}
\end{lemma}

\begin{proof}
The paired third-cumulant identity \eqref{eq:paired-cumulants} gives
$V_A(\lambda)\ge V_A(2/N)$. Lemma~\ref{lem:u2-continuum-variance}
at $M=N/2$ now gives
\[
 V_A(\lambda)>\frac{N^3}{8}\left(\frac{57}{100}-\frac6N\right)
 \ge\frac{69}{1000}N^3 .
\]
For the upper bound, \eqref{eq:variance-origin} gives
$V_A(\lambda)\le(2N^3+N^2)/12<N^3/5$.
Denote the canceled numerator and denominator sets by
$\mathcal U$ and $\mathcal V$, respectively.

Every Euler index is at most $2N$, so $A(e^{-z})$ has no zero in the
strip $|\Im z|<\pi/N$. At $z=0$ the common Euler zeros cancel and
$A(1)>0$. Take the analytic logarithm real on the real axis and put
$\mathcal L(z)=\log A(e^{-z})+Dz/2$. Reciprocity makes $\mathcal L$ even, and the
centered logarithm satisfies
\[
 \log\phi(\theta)=\mathcal L(\lambda-i\theta)-\mathcal L(\lambda)+i\theta \mathcal L'(\lambda).
\]
By \eqref{eq:sinh-product},
\[
 \frac{\sinh(jz/2)}{jz/2}
 =\prod_{\ell\ge1}\left(1+\frac{j^2z^2}{4\pi^2\ell^2}\right),
\]
the logarithmic derivative of each $(j,\ell)$ factor has poles at
$\pm2\pi i\ell/j$. Their distances from any real $s\ge0$ are at least
$2\pi\ell/j$. Absolute summation of the differentiated series gives,
for $r\ge2$,
\[
 |\mathcal L^{(r)}(s)|\le
 \frac{2(r-1)!\zeta(r)}{(2\pi)^r}
 \sum_{j\in\mathcal U\cup\mathcal V}j^r .
\]
For odd $r\ge3$, $\mathcal L^{(r)}(0)=0$. Integrating the bound of order
$r+1$ on $[0,\lambda]$ gives
\[
 |\mathcal L^{(r)}(\lambda)|
 \le\frac{2r!\zeta(r+1)}{(2\pi)^{r+1}}\,
 \lambda\sum_{j\in\mathcal U\cup\mathcal V}j^{r+1}.
\]

Put $S_4=\sum_{j\in\mathcal U\cup\mathcal V}j^4$.
The Taylor series at $\lambda$ has radius at least $\pi/N$.
The even terms starting at order four contribute to $R$ and the odd
terms starting at order three to $I$. Using
$\sum j^r\le(2N)^{r-4}S_4$ and $\zeta(r)\le\zeta(4)$ for $r\ge4$,
and summing the resulting geometric series, we find
\[
 |R|\le\frac{S_4\theta^4}{2880(1-(N\theta/\pi)^2)},\qquad
 |I|\le\frac{\lambda S_4|\theta|^3}{720(1-(N\theta/\pi)^2)} .
\]
The two index sets are disjoint, whence
\[
 S_4\le\sum_{j=1}^{2N}j^4
 =\frac{32}{5}N^5+8N^4+\frac83N^3-\frac1{15}N
 <\frac{321}{50}N^5.
\]
For $N|\theta|\le h$ and $\pi^2>49/5$ the geometric denominator
exceeds $11/16$. This proves the bounds on $R$ and $I$.
The last bound in \eqref{eq:small-local} follows from
$1/2-(13/4000)h^2/(69/1000)>71/200$.

Write $v=69/1000$, $a=71/200$, $c_4=13/4000$, and $c_3=13/1000$.
Applying $e^R\ge1-|R|$ and $1-\cos I\le I^2/2$ in
$\Re\phi=e^{\Re\log\phi}-e^{\Re\log\phi}(1-\cos I)$ gives
\[
 \Re\phi\ge
 e^{-V\theta^2/2}(1-c_4N^5\theta^4)
 -\tfrac12 c_3^2\lambda^2N^{10}\theta^6 e^{-aV\theta^2}.
\]
On the local arc $0\le\cos(\theta/2)\le1$ and
$\cos(\theta/2)\ge1-\theta^2/8$. By \eqref{eq:half-moments} the fourth-order and phase losses,
divided by $I_G$, are respectively at most
\[
 \frac{3c_4}{v^2N}<\frac{2.048}{N},\qquad
 \frac{15c_3^2}{4a^3v^3\sqrt{2a}\,N}<\frac{51.175}{N},
\]
where $V\ge vN^3$ and $\lambda\le2/N$.
The cosine loss is at most $1/(8V)<1/(100N)$. For the omitted
Gaussian tail let $z=h\sqrt V/N>1$. Its relative mass is at most
$\sqrt{2/\pi}\,e^{-z^2/2}/z<e^{-vh^2N/2}$. The function
$Ne^{-vh^2N/2}$ decreases for $N\ge500$, and at $500$ it is less
than $500/2^{52}<1/10$. The four losses sum to less than $54/N$.
For the imaginary part,
$|\Im\phi\sin(\theta/2)|\le e^{-aV\theta^2}|I|\theta/2$.
The fourth moment gives the relative bound
\[
 \frac{3c_3}{8a^2v^2\sqrt{2a}}\frac{\lambda}{N}
 \le\frac{10\lambda}{N}.
\]
This is strict for $\lambda>0$, and $I=0$ for $\lambda=0$.

\end{proof}

At the saddle $\mu_A(\lambda)=k-1/2$, the Fourier factors of the
two adjacent coefficients differ by
\[
 \Delta_\lambda(\theta)
 = e^{\lambda/2-i\theta/2}-e^{-\lambda/2+i\theta/2}
 =2\sinh(\lambda/2)\cos(\theta/2)
  -2i\cosh(\lambda/2)\sin(\theta/2).
\]
In fact, coefficient inversion and conjugate symmetry give
\begin{align*}
 a_k-a_{k-1}
 &=\frac{A(q)e^{\lambda(k-1/2)}}{2\pi}
   \int_{-\pi}^{\pi}\phi(\theta)\Delta_\lambda(\theta)\,d\theta\\
 &=\frac{A(q)e^{\lambda(k-1/2)}}{\pi}
   \int_0^\pi K_\lambda(\theta)\,d\theta,\qquad
 K_\lambda(\theta)=\Re\bigl(\phi(\theta)\Delta_\lambda(\theta)\bigr).
\end{align*}
Thus the two integrals just estimated occur with respective
multipliers $2\sinh(\lambda/2)$ and $2\cosh(\lambda/2)$.

\begin{lemma}\label{lem:paired-first-lobe}
Suppose $N\ge500$, $1\le k\le\lfloor D/2\rfloor$, and the solution of
$\mu_A(\lambda)=k-1/2$ satisfies $\lambda N\le2$. Then
\begin{equation}\label{eq:positive-first-lobe}
 \int_0^{\pi/N}K_\lambda(\theta)\,d\theta
 >\frac35\lambda I_G .
\end{equation}
In particular $a_k>a_{k-1}$ whenever
\begin{equation}\label{eq:far-criterion}
 \int_{\pi/N}^{\pi}(\lambda+\theta)|\phi(\theta)|\,d\theta
 \le\frac12\lambda I_G .
\end{equation}
\end{lemma}

\begin{proof}
Since $\mu_A(0)=D/2$ and $-\mu_A'=V_A$,
\[
 \frac12\le D/2-k+1/2=\int_0^\lambda V_A(s)\,ds
 \le\frac15\lambda N^3.
\]
Thus $\lambda\ge5/(2N^3)$ and $I_G>(5/2)N^{-3/2}$.
Use $2\sinh(\lambda/2)\ge\lambda$, $\cosh(\lambda/2)<9/8$, and
\eqref{eq:local-cos}--\eqref{eq:local-sin} to obtain
\[
 \int_0^{h/N}K_\lambda(\theta)\,d\theta
 \ge \left(1-\frac{54+2(9/8)10}{N}\right)\lambda I_G
 >\left(1-\frac{77}{N}\right)\lambda I_G .
\]

For $h/N\le\theta\le\pi/N$, the parameters of
Lemma~\ref{lem:first-lobe-separated-pairs} satisfy
\[
 a\ge\left(\frac7{16}-\frac16\left(\frac7{16}\right)^3\right)^2,
 \qquad b\le\sinh^2(1/2)<\frac{11}{40}.
\]
The quantity $a(3+2b-2a)/(3+4b)^2$ is nondecreasing in $a$ and
nonincreasing in $b$ for $0\le a\le1/2$, $b\ge0$. The displayed
endpoints give a lower bound greater than $17/500$.
Since $4\lfloor N/4\rfloor\ge N-4$,
$|\phi(\theta)|\le e^{-17(N-4)/500}$ on this arc.
Also, for $0\le\theta\le\pi$,
\[
 |K_\lambda(\theta)|\le|\phi(\theta)|\,|\Delta_\lambda(\theta)|
 \le\tfrac98(\lambda+\theta)|\phi(\theta)| ,
\]
by $\sinh x\le x\cosh x$ and $|\sin(\theta/2)|\le\theta/2$.
Integrating and applying the lower bounds on $\lambda$ and $I_G$ yields
\[
 \frac{1}{\lambda I_G}\int_{h/N}^{\pi/N}|K_\lambda(\theta)|\,d\theta
 \le\frac98 e^{-17(N-4)/500}
 \left(\frac{2(\pi-h)}5N^{1/2}
       +\frac{2(\pi^2-h^2)}{25}N^{5/2}\right)<\frac15.
\]
Both summands decrease for $N\ge500$. At $500$ the last comparison
uses $\pi<22/7$, $\sqrt{500}<23$, and $e^{2108/125}>20{,}000{,}000$
which follows by truncating the exponential series at order $45$.
Now $1-77/500-1/5>3/5$ gives \eqref{eq:positive-first-lobe}.
Under \eqref{eq:far-criterion} the rest of the integral has absolute
value at most $(9/16)\lambda I_G<(3/5)\lambda I_G$. The inversion
formula then gives $a_k>a_{k-1}$.
\end{proof}

The preceding estimates now give positivity for the entire range
$0<N\lambda\le2$.

\begin{lemma}\label{lem:paired-saddle-positive}
Let $N\ge1000$, $0\le m\le N$, and
$1\le k\le\lfloor D/2\rfloor$. If the unique positive solution of
\[
 \mu_A(\lambda)=k-\frac12
\]
satisfies $\lambda N\le2$, then $a_k>a_{k-1}$.
\end{lemma}
\begin{proof}
Use the notation of Lemma~\ref{lem:paired-first-lobe}.
For a variable tilt $s$ set
\[
 \phi_s(\theta)=e^{-i\mu_A(s)\theta}
     \frac{A_{m,N}(e^{-s+i\theta})}{A_{m,N}(e^{-s})}.
\]
Then $\phi=\phi_\lambda$. First suppose $1/4\le\lambda N\le2$.
On $\pi/N\le\theta\le100/N$, we have $3<N\theta\le100$ and
$\int_{\pi/N}^{100/N}(\lambda+\theta)\,d\theta
\le(2\cdot100+100^2/2)/N^2=5200/N^2$.
Lemma~\ref{lem:compact-arc} thus gives
\[
 \int_{\pi/N}^{100/N}(\lambda+\theta)|\phi(\theta)|\,d\theta
 \le5200N^{-2}e^{-N/16}<\frac{\lambda I_G}{100}.
\]
Indeed, $I_G>(5/2)N^{-3/2}$ and $\lambda\ge1/(4N)$, hence
$\lambda I_G>(5/8)N^{-5/2}$. It suffices to check
$8320\sqrt N e^{-N/16}<1/100$. The left side decreases for
$N\ge1000$. At $N=1000$ use $\sqrt{1000}<32$ and $e^{125/2}>2^{62}$.
On $100/N\le\theta\le\pi$, one has $\lambda+\theta<4$.

On $100/N\le\theta\le\pi$ we have
$\csc(\theta/2)\le\csc(50/N)\le48N/2399$, by
$\sin x\ge x-x^3/6$ and $50/N\le1/20$.
Since $26/625-117/11995>159/5000$,
Lemma~\ref{lem:wz-outer-bound} gives
\[
 2\int_{100/N}^{\pi}|\phi(\theta)|\,d\theta
 \le208\pi\sqrt{N+1}\,e^{-159N/5000}.
\]
The bounds $V\le V_A(0)\le N^3/4$ and $\lambda\ge1/(4N)$ give
$\lambda I_G\ge(\sqrt{2\pi}/4)N^{-5/2}$. The ratio is therefore at most
\[
 416\sqrt{2\pi}\,N^3\sqrt{1+1/N}\,e^{-159N/5000}.
\]
Its logarithmic derivative is less than $3/N-159/5000<0$.
At $N=1000$, use $\sqrt{2\pi}<251/100$, $\sqrt{1+1/1000}<1001/1000$,
and $\sum_{j=0}^{40}(159/5)^j/j!>6\cdot10^{13}$. The resulting
bound is $416(251/100)(1001/1000)/60000=6532526/375000000<1/50$.

Since $\lambda+\theta<4$, it follows that
\[
 \int_{100/N}^{\pi}(\lambda+\theta)|\phi(\theta)|\,d\theta
 <\frac{\lambda I_G}{25}.
\]
The sum is less than $\lambda I_G/2$, as required.

Now suppose $0<\lambda N\le1/4$. Apply
Lemma~\ref{lem:strip-transfer} with
$\varepsilon=1/4$ and $\theta_0=\pi/N$. At the vertical boundary
$s=1/(4N)$ the compact arc is bounded by $e^{-N/16}$.
On $\theta\ge100/N$, Lemma~\ref{lem:wz-outer-bound} and
$\csc(\theta/2)\le48N/2399$ give the same bound, because
\[
 \frac{12}{125}-\frac{39}{80}\frac{48}{2399}-\frac1{16}
 >\frac{23}{1000},\qquad
 64\sqrt{N+1}e^{-23N/1000}<1\quad(N\ge1000).
\]
For the last inequality, the left side decreases with $N$, and its
value at 1000 is less than $2048e^{-23}<1$.
Conjugation supplies the rest of the vertical boundary.

At $\theta=\pi/N$, Lemma~\ref{lem:first-lobe-separated-pairs}
has $a=1/2$ and $b=\sinh^2(sN/4)<1/100$ for $0\le sN\le1/4$,
since the power series gives $\sinh(1/16)\le\sum_{j\ge1}(1/16)^j=1/15$.
Thus
\[
 \frac{a(3+2b-2a)}{(3+4b)^2}
 =\frac{1+b}{(3+4b)^2}\ge\frac1{(3+4/100)^2}>\frac1{10}.
\]
Since $4\lfloor N/4\rfloor\ge N-4$, the boundary bound is
$|\phi_s(\pi/N)|\le e^{-(N-4)/10}\le e^{-N/16}$. The transfer factor
is at most $\exp(N/160)$. The strip transfer consequently gives the uniform bound
\[
 |\phi_s(\theta)|\le e^{-9N/160}
 \quad\left(0\le sN\le\frac14,\quad\frac\pi N\le\theta\le\pi\right).
\]
Since $\pi\lambda+\pi^2/2<6$, integration gives
\[
 \int_{\pi/N}^\pi(\lambda+\theta)|\phi(\theta)|\,d\theta
 <6e^{-9N/160}<\frac{\lambda I_G}{2}.
\]
Here $\lambda I_G>(25/4)N^{-9/2}$ by the half-integer saddle
lower bound. The function $N^{9/2}e^{-9N/160}$ decreases for
$N\ge1000$. At $N=1000$ the last comparison follows from
$\sqrt{1000}<32$ and $e^{225/4}>2^{56}$.
Thus the far criterion holds in both tilt ranges, proving positivity.
\end{proof}

\section{The saddle for the first-difference polynomial}\label{sec:large}
We now use $k=\mu_B(\lambda)$. The range $N\lambda\ge4$ admits
one argument with two sets of constants. We then extend positivity
to $N\lambda\ge199/100$, overlapping Section~\ref{sec:small}.

\subsection{A common estimate for large tilts}
For a lower bound we must account for the interval $(m,2m]$ missing
from the partition product. Its largest possible loss is controlled
by the following integral.

\begin{lemma}\label{lem:missing-interval}
For $f_2(x)=(x/(2\sinh(x/2)))^2$, with $f_2(0)=1$, one has
\[
 \max_{s\ge0}\int_s^{2s}f_2(x)\,dx<\frac{51}{50}.
\]
The unique positive maximizing point is $s_*=2\log(1+\sqrt2)$.
\end{lemma}
\begin{proof}
Since $f_2(2s)=f_2(s)/\cosh^2(s/2)$,
\[
 \frac{d}{ds}\int_s^{2s}f_2(x)\,dx
 =f_2(s)\left(\frac{2}{\cosh^2(s/2)}-1\right).
\]
The derivative is positive for $0<s<s_*$ and negative for $s>s_*$.
Put
\[
 s_-:=\frac{17627}{10000},\qquad s_+:=\frac{4407}{2500},\qquad
 r_-:=\frac{4289}{25000},\qquad r_+:=\frac{429}{2500}.
\]
The following rational bounds locate $s_*$:
\[
 \begin{aligned}
 \cosh(s_-/2)
 &\le\sum_{j=0}^{5}\frac{(s_-/2)^{2j}}{(2j)!}
   +\frac{182}{181}\frac{(s_-/2)^{12}}{12!}
   <\frac{7071}{5000}<\sqrt2,\\
 \cosh(s_+/2)
 &\ge\sum_{j=0}^{4}\frac{(s_+/2)^{2j}}{(2j)!}
   >\frac{70711}{50000}>\sqrt2.
 \end{aligned}
\]
Here the factor $182/181$ bounds the tail because $s_-/2<1$ and,
after the term of degree $12$, the ratio of consecutive terms is
at most $1/(13\cdot14)$. The two bounds on $\sqrt2$ also give
\[
 s_-<s_*<s_+,\qquad r_-<e^{-s_*}=3-2\sqrt2<r_+.
\]
The nonnegative expansion
$f_2(x)=\sum_{\ell\ge1}\ell x^2e^{-\ell x}$ permits termwise
integration. For $s=s_*$ and $r=e^{-s_*}$ it yields
\[
 \int_s^{2s}f_2(x)\,dx
 =\sum_{\ell\ge1}
 \left\{r^\ell\left(s^2+\frac{2s}{\ell}+\frac{2}{\ell^2}\right)
 -r^{2\ell}\left(4s^2+\frac{4s}{\ell}+\frac{2}{\ell^2}\right)\right\}.
\]
For the first three terms, the endpoint bounds give
\[
 \begin{aligned}
 r_+(s_+^2+2s_++2)-r_-^2(4s_-^2+4s_-+2)
 &<\frac{84924}{100000},\\
 r_+^2(s_+^2+s_++1/2)-r_-^4(4s_-^2+2s_-+1/2)
 &<\frac{14389}{100000},\\
 r_+^3(s_+^2+2s_+/3+2/9)-r_-^6(4s_-^2+4s_-/3+2/9)
 &<\frac{2239}{100000}.
 \end{aligned}
\]
Each displayed left side is the upper bound obtained by evaluating
the positive term at $(s_+,r_+)$ and the subtracted term at
$(s_-,r_-)$.
For the remaining terms, discard the subtracted terms and use
$\ell\ge4$ to obtain
\[
 \sum_{\ell\ge4}r^\ell
 \left(s^2+\frac{2s}{\ell}+\frac{2}{\ell^2}\right)
 \le\frac{r_+^4}{1-r_+}
       \left(s_+^2+\frac{s_+}{2}+\frac18\right)
 <\frac{431}{100000}.
\]
Consequently
\[
 \int_{s_*}^{2s_*}f_2(x)\,dx
 <\frac{84924+14389+2239+431}{100000}
 =\frac{101983}{100000}<\frac{51}{50}.
\]
\end{proof}

The remaining variance loss comes from indices greater than $N$.
On this tail the ratio of the third to the second cumulant is large
enough to bound both cumulants at once.

\begin{lemma}\label{lem:large-tilt-sharp-variance}
Let $0\le m\le n=N$, $\lambda>0$, $M=1/\lambda\ge17$, and
$u=\lambda N\ge8$. Define
\[
 B(q)=(1-q)A_{m,n}(q),\qquad
 V_B(\lambda)=\frac{d^2}{d\lambda^2}\log B(e^{-\lambda}),
 \qquad K_{3,B}(\lambda)=-\frac{d^3}{d\lambda^3}
                                      \log B(e^{-\lambda}).
\]
The real logarithm is well defined because $B(q)>0$ for $0<q<1$.
Then
\begin{align*}
 V_B(\lambda)&\ge
 \left(\frac{\pi^2}{3}-1.02-\frac2M-0.055054\right)M^3
 >2.0971M^3,\\
 K_{3,B}(\lambda)&\le3.584M V_B(\lambda).
\end{align*}
\end{lemma}
\begin{proof}
The variance contribution of an inverse Euler factor $(1-e^{-j\lambda})^{-1}$
is $M^2f_2(j/M)$, and its third cumulant is $M^3f_3(j/M)$.
Let $E_0=\{j\ge1:j\ne1,\ j\notin(m,2m]\}$ and
$L_0(\lambda)=-\sum_{j\in E_0}\log(1-e^{-j\lambda})$.
Put $V_0=L_0''$ and $K_{3,0}=-L_0'''$. Let
\[
 z=\frac1M\left(f_2(1/M)+\sum_{m<j\le2m}f_2(j/M)\right).
\]
If $m=0$, the removed interval is empty. In all cases monotonicity,
Lemma~\ref{lem:missing-interval}, and $f_2\le1$ give
$0\le z\le1.02+1/M$. For a positive decreasing function $f$ the right-endpoint sum satisfies
\[
 \int_0^\infty f(x)\,dx-\frac{f(0)}M
 \le \frac1M\sum_{j\ge1}f(j/M)\le\int_0^\infty f(x)\,dx .
\]
Apply this inequality to $f_2,f_3$, using $f_3\ge2f_2$ on the removed
indices. It follows that
\[
 V_0\ge M^3(\frac{\pi^2}{3}-z-1/M),\qquad
 K_{3,0}\le M^4(\pi^2-2z).
\]
Since the ratio $(\pi^2-2z)/(\frac{\pi^2}{3}-z-1/M)$ increases with $z$, it follows that
\[
 K_{3,0}\le3.584M V_0\qquad(M\ge17).
\]
Indeed the denominator is at least $\frac{\pi^2}{3}-1.02-2/17>0$, and the ratio
also increases in $1/M$ with $z$ fixed. Its value at
$M=17$, $z=1.02+1/17$, is less than $3.584$.

In the canceled representation of $B$, its inverse factors are
indexed by $E_0\cap[1,N]$ and its numerator factors by $\mathcal U\subset(N,2N]$.
To obtain its signed cumulants from those of $L_0$, subtract first the
inverse-factor contributions of $E_0\cap(N,\infty)$ and then those
of $\mathcal U$. An index in both sets must be subtracted twice. Thus the total
subtraction is bounded by two full tails $j>N$. Since $f_2$ is decreasing,
\[
 \sum_{j>N}M^2f_2(j/M)
 \le M^3\int_u^\infty f_2(x)\,dx
 \le M^3\frac{e^{-u}(u^2+2u+2)}{(1-e^{-u})^2}
 <0.027527M^3.
\]
This proves the lower bound on $V_B$. Moreover, at every such removed
index,
\[
 \frac{M^3f_3(j/M)}{M^2f_2(j/M)}
 =j\coth(j/(2M))>8M>3.584M.
\]
Thus removing these contributions can only decrease
$K_3-3.584MV$. The already proved inequality for $V_0,K_{3,0}$ therefore
proves the third-cumulant comparison for $B$.
\end{proof}

To bound the product away from its saddle, compare its modulus with
the partition generating function
\[
 P(q)=\prod_{j\ge1}(1-q^j)^{-1},\qquad |q|<1 .
\]
The positivity of the central Gaussian polynomial controls the same
missing interval that arose in the variance estimate.

\begin{lemma}\label{lem:half-partition-comparison}
Let $0\le m\le n=N$ with $N\ge1$, let $q=e^{-\lambda}$ with $\lambda>0$, and put
$M=\lambda^{-1}$ and $u=\lambda N$. Define
\begin{equation}\label{eq:partition-correction}
 \epsilon_*(u)=3\sum_{\substack{\ell\ge1\\\ell\ {\rm odd}}}
                    \frac{e^{-\ell u}}{\ell^2}.
\end{equation}
For every real $\theta$ we have
\begin{equation}\label{eq:half-partition-bound}
 \left|\frac{A_{m,n}(qe^{i\theta})}{A_{m,n}(q)}\right|
 \le e^{M\epsilon_*(u)}
 \left|\frac{P(qe^{i\theta})}{P(q)}\right|^{1/2}.
\end{equation}
Consequently,
\[
 \left|\frac{A_{m,n}(qe^{i\theta})}{A_{m,n}(q)}\right|
 \le
 \exp\!\left(\frac{3M e^{-u}}{1-e^{-u}}\right)
 \left|\frac{P(qe^{i\theta})}{P(q)}\right|^{1/2}.
\]
In particular, if $u\ge8$, the exponential correction is at most
$\exp(0.001007M)$.
\end{lemma}
\begin{proof}
For $j\ge1$, define
\[
 \ell_j=\log\frac{|1-q^je^{ij\theta}|}{1-q^j}\ge0,
 \qquad \mathcal S_N=\sum_{j\le N}\ell_j,
 \qquad \mathcal S_\infty=\sum_{j\ge1}\ell_j.
\]
The bound $\ell_j\le\log((1+q^j)/(1-q^j))=O(q^j)$ proves convergence.
The nonnegative coefficients of
$\qbinom{2m}{m}_q$ give
\[
 \left|\frac{\qbinom{2m}{m}_{qe^{i\theta}}}{\qbinom{2m}{m}_q}\right|\le1.
\]
Taking logarithms in its product formula yields the exact inequality
\[
 \sum_{m<j\le2m}\ell_j\le\sum_{j\le m}\ell_j.
\]
In both cases of the canceled-factor table, the denominator set is
$\mathcal V=[1,N]\setminus(m,2m]$ and $[1,m]\subseteq\mathcal V$. Hence
\[
 \sum_{j\in[1,N]\setminus\mathcal V}\ell_j\le\sum_{m<j\le2m}\ell_j
 \le\sum_{j\le m}\ell_j\le\sum_{j\in\mathcal V}\ell_j,
 \qquad \sum_{j\in\mathcal V}\ell_j\ge \mathcal S_N/2.
\]
All numerator indices are greater than $N$. Therefore
\[
 \log\left|\frac{A_{m,n}(qe^{i\theta})}{A_{m,n}(q)}\right|
 =\sum_{j\in\mathcal U}\ell_j-\sum_{j\in\mathcal V}\ell_j
 \le-\frac12 \mathcal S_\infty+\frac32\sum_{j>N}\ell_j .
\]
Set $g(x)=\log((1+e^{-x})/(1-e^{-x}))$. Since $\ell_j\le g(j/M)$
and $g$ is decreasing,
\[
 \sum_{j>N} \ell_j\le M\int_u^\infty g(x)\,dx
 =2M\sum_{\substack{\ell\ge1\\ \ell\ {\rm odd}}}
                     \frac{e^{-\ell u}}{\ell^2}.
\]
The equality follows by integrating the nonnegative power series term
by term. Also $-\mathcal S_\infty=\log|P(qe^{i\theta})/P(q)|$ and
$\epsilon_*(u)\le3\sum_{\ell\ge1}e^{-\ell u}=3/(e^u-1)$.
This gives both inequalities. For $u\ge8$, use $3/(e^8-1)<0.001007$.
\end{proof}

For the partition product, summation of two geometric series gives
for $a>0$
\[
 S(a,t)=\sum_{j\ge1}e^{-aj}(1-\cos jt)
 =\frac12\coth(a/2)\frac{1-\cos t}{\cosh a-\cos t}.
\]
Expanding the logarithm of $P$ gives
\[
 \mathcal F(M,\theta):=-\log\left|\frac{P(e^{-1/M+i\theta})}{P(e^{-1/M})}\right|
 =\sum_{\ell\ge1}\frac1\ell S(\ell/M,\ell\theta).
\]
Write
\[
 \omega(M,\theta)=\tfrac12\log\left(1+
        \frac{\sin^2(\theta/2)}{\sinh^2(1/(2M))}\right),\qquad
 \mathcal E_{M,\epsilon}(\theta)=
        \exp\{\omega(M,\theta)-\mathcal F(M,\theta)/2+\epsilon M\}.
\]
The factor $e^{\omega}$ is $|1-qe^{i\theta}|/(1-q)$ at $q=e^{-1/M}$.
Thus Lemma~\ref{lem:half-partition-comparison} gives
$|B(qe^{i\theta})/B(q)|\le \mathcal E_{M,\epsilon}(\theta)$ whenever
$3/(e^u-1)\le\epsilon$. On the ranges below we use
\[
 e^{\omega(M,\theta)}\le\sqrt{1+(M\theta)^2}
 \quad\hbox{and}\quad e^{\omega(M,\theta)}\le1+2M\sin(\theta/2).
\]
The first bound follows from $\sin(\theta/2)\le\theta/2$ and
$\sinh(1/(2M))\ge1/(2M)$; the second follows from
$\sqrt{1+x^2}\le1+x$ for $x\ge0$.

\smallskip
We next bound the local integral uniformly in the Euler index sets.
Put $\mathcal V_0=\mathcal V\setminus\{1\}$ and $V=V_B(\lambda)$. Define
\[
 S_{r,E}=\sum_{j\in E}\sum_{\ell\ge1}j^r\ell^{r-1}e^{-j\ell/M}.
\]
Equivalently, $S_{r,E}$ is $(-1)^r$ times the $r$th derivative of
$-\sum_{j\in E}\log(1-e^{-j\lambda})$. In particular
$V=S_{2,\mathcal V_0}-S_{2,\mathcal U}$ and $K_{3,B}=S_{3,\mathcal V_0}-S_{3,\mathcal U}$.
The kernel for $r=4$ is
\[
 f_4(x)=x^4\sum_{\ell\ge1}\ell^3e^{-\ell x}
       =x^4\frac{e^{-x}(1+4e^{-x}+e^{-2x})}{(1-e^{-x})^4},\qquad
 \int_0^\infty f_4(x)\,dx=24\zeta(2)=4\pi^2 .
\]
We need an upper bound for its full discrete sum and its tail.

To determine the shape of $f_4$, put $y=x/2$, $s=\sinh y$, and $c_y=\cosh y$.
The sign of $f_4'$ is the opposite of the sign of
$h(y)=yc_y(3+s^2)-s(3+2s^2)$. The derivative chain is
\[
 h'=s g,\quad g=y(5+3s^2)-5sc_y,\quad
 g'=s\ell,\quad \ell=6yc_y-7s,\quad
 \ell'=6ys-c_y,\quad \ell''=5s+6yc_y>0.
\]
Since $\ell'(0)=-1$ and $\ell''(y)>0$ for $y>0$, $\ell'$ crosses
zero once. Each of $\ell,g,h$ vanishes at zero and tends to infinity.
The derivative identities now show successively that each crosses from
negative to positive exactly once on the positive axis. Thus $f_4$ is
unimodal.
Also, with $z=y^2/6$ and $\sinh y\ge y+y^3/6$,
\[
 f_4=6(y/\sinh y)^4+4y^4/\sinh^2y
 \le\frac6{(1+z)^4}+\frac{24z}{(1+z)^2}<7.
\]
The last inequality follows because the numerator of the difference is
$1+2z(2z^2-3z+2)+7z^4>0$.
The inequality $f_4\ge6f_2$ follows from
$f_4/f_2=x^2+6f_2$ and
$f_2(x)\ge1-x^2/12$. For the latter bound, the coefficient of $y^2$
in $(1-y^2/3)\sinh^2 y$ is one and the coefficient of $y^{2k}$,
$k\ge2$, has the sign of $4/((2k)(2k-1))-1/3\le0$.
For $0<t<7$ the superlevel set $\{x\ge0:f_4(x)>t\}$ is an
interval. If its length is $d_t$, it contains at most $Md_t+1$
points of the grid $\{j/M:j\ge1\}$. Integrating this count over $t$
gives
\[
 \sum_{j\ge1}f_4(j/M)\le M\int_0^\infty f_4(x)\,dx+7=4\pi^2M+7.
\]

For $r=e^{-x}$,
\[
 (\log f_4)'(x)=\frac4x-1-\frac{4r+2r^2}{1+4r+r^2}
                         -\frac{4r}{1-r}<\frac4x-1 .
\]
Consequently all three kernels $f_2,f_3,f_4$ are decreasing for $x\ge4$.
For $u=N/M\ge4$ put $r_u=e^{-u}$ and
$T_d=M^d\sum_{j>N} f_d(j/M)$. Integrating $x^d e^{-x}$ from $u$ to
infinity and bounding the remaining factor at $u$ yields
\begin{align*}
 T_2/M^3&\le\frac{r_u}{(1-r_u)^2}(u^2+2u+2),\\
 T_3/M^4&\le\frac{r_u(1+r_u)}{(1-r_u)^3}(u^3+3u^2+6u+6),\\
 T_4/M^5&\le\frac{r_u(1+4r_u+r_u^2)}{(1-r_u)^4}
             (u^4+4u^3+12u^2+24u+24).
\end{align*}
These bounds decrease with $u$: the product of $r_u$ with each
polynomial is the tail integral and the factor left over also decreases.
Return to the missing-index sum
$z=M^{-1}(f_2(1/M)+\sum_{m<j\le2m} f_2(j/M))$
of Lemma~\ref{lem:large-tilt-sharp-variance}. The denominator
set $\mathcal V_0$ is contained in $E_0=\{j\ge1:j\ne1,\ j\notin(m,2m]\}$.
Removing the missing indices from the full fourth-cumulant sum and
using $f_4\ge6f_2$ gives the first of the bounds
\[
 S_{4,\mathcal V_0}\le M^5(4\pi^2-6z+7/M),\qquad
 V\ge M^3(\frac{\pi^2}{3}-z-1/M)-2T_2,\qquad 0\le z\le1.02+1/M .
\]
This second inequality follows by taking out the two variance tails in
the same lemma. Also $S_{d,\mathcal U}\le T_d$ for $d=2,3,4$.

At $k=\mu_B(\lambda)$ put
\[
 \phi(\theta)=e^{-ik\theta}\frac{B(e^{-\lambda+i\theta})}{B(e^{-\lambda})},
 \qquad \log\phi=R+iI,\qquad R(0)=I(0)=0.
\]
There is a continuous logarithm because all Euler factors are nonzero
on $|q|=e^{-\lambda}$. Write $\sigma_j=1$ on $\mathcal V_0$ and $-1$ on $\mathcal U$.
The saddle identity cancels the linear term in
\[
 R+iI=\sum_{j\in \mathcal V_0\cup\mathcal U}\sigma_j\sum_{\ell\ge1}
        \frac{e^{-j\ell/M}}{\ell}
        \big(e^{ij\ell\theta}-1-ij\ell\theta\big).
\]
Applying $x^2/2-x^4/24\le1-\cos x\le x^2/2$ and
$|\sin x-x|\le |x|^3/6$ gives
\begin{align*}
 -V\theta^2/2-S_{4,\mathcal U}\theta^4/24\ &\le R
      \le -V\theta^2/2+S_{4,\mathcal V_0}\theta^4/24,\\
 |I|\ &\le (S_{3,\mathcal V_0}+S_{3,\mathcal U})|\theta|^3/6 .
\end{align*}
Suppose $V>0$ and that $R\le-aV\theta^2$ and
$|I|\le sMV\theta^3/6$ on $0\le\theta\le\theta_0$, where $a,s>0$.
From $\Re\phi=e^R-e^R(1-\cos I)$, $1-\cos I\le I^2/2$ and
$e^{-x}\ge1-x$, we obtain
\[
 \Re\phi(\theta)\ge
 e^{-V\theta^2/2}\bigl(1-S_{4,\mathcal U}\theta^4/24\bigr)
 -\frac{s^2M^2V^2\theta^6}{72}e^{-aV\theta^2}.
\]
Put $\eta(z)=\sqrt{2/\pi}\,e^{-z^2/2}/z$ for $z>0$. Integrating
the two errors to infinity gives
\begin{equation}\label{eq:shared-gaussian-loss}
 \frac1{I_G}\int_0^{\theta_0}\Re\phi(\theta)\,d\theta
 \ge 1-\eta(\theta_0\sqrt V)-\frac{S_{4,\mathcal U}}{8V^2}
       -\frac{15s^2M^2}{576a^3\sqrt{2a}\,V}.
\end{equation}
Indeed, \eqref{eq:half-moments} gives
$J_2(1/2)/I_G=3/V^2$ and $J_3(a)/I_G=15/(8a^3\sqrt{2a}\,V^3)$.
For $t\ge z$ bound $e^{-t^2/2}$ by $(t/z)e^{-t^2/2}$ to obtain
$\eta(\theta_0\sqrt V)$ for the omitted relative Gaussian tail.

\begin{lemma}\label{lem:partition-tail-elementary}
For $M\ge17$ and $0\le\epsilon\le0.001007$,
\[
 \int_{1/(2M)}^\pi \mathcal E_{M,\epsilon}(\theta)\,d\theta<0.55I_*(M).
\]
For $M\ge125$ and $0\le\epsilon\le0.055973$,
\[
 \int_{2/(5M)}^\pi \mathcal E_{M,\epsilon}(\theta)\,d\theta<0.006I_*(M).
\]
\end{lemma}
\begin{proof}
We first treat $c/M\le\theta\le1/M$, where $c=1/2$ or $2/5$.
For $1\le j\le17$, both $x=j/M$ and $y=j\theta$ lie in $[0,1]$.
The inequalities
\[
 \coth(x/2)\ge2/x,\qquad
 \cosh x-1\le\frac{11}{20}x^2,\qquad
 1-\cos y\ge\frac{11}{24}y^2
\]
give
\[
 S(j/M,j\theta)\ge\frac{M}{j}
       \frac{(M\theta)^2}{6/5+(M\theta)^2}.
\]
The positive Taylor coefficients show that
$(\cosh x-1)/x^2$ increases for $x\ge0$. Thus the hyperbolic bound follows from
$\cosh1-1<1/2+1/24+(1/720)/(1-1/56)<11/20$.
Also $\sum_{j=1}^{17}j^{-2}>\pi^2/6-1/17>79/50$.
With $\tau=M\theta$ and $h(\tau)=79\tau^2/(120+100\tau^2)$,
we therefore have
\[
 \frac12\mathcal F(M,\theta)\ge Mh(\tau),\qquad
 e^{\omega(M,\theta)}\le\sqrt{1+\tau^2}.
\]
Since
\[
 h'(\tau)=\frac{(237/125)\tau}{(6/5+\tau^2)^2},
\]
its unique positive maximum occurs at $\sqrt{2/5}$, and its
values at $2/5$ and $1$ both exceed $3/8$.
Thus $h(\tau)\ge h(c)+(3/8)(\tau-c)$ for $c\le\tau\le1$, so
\begin{equation}\label{eq:elementary-short-partition-tail}
 \frac1{I_*(M)}\int_{c/M}^{1/M}\mathcal E_{M,\epsilon}(\theta)\,d\theta
 \le\frac{8\sqrt2}{3\sqrt{3/(2\pi)}\sqrt M}
       e^{-(h(c)-\epsilon)M}.
\end{equation}
This bound decreases with $M$ in each of the two ranges of the
lemma. In the first range
$17(h(1/2)-0.001007)>229/100$ and
$e^{229/100}>\sum_{j=0}^{12}(229/100)^j/j!>49/5$.
Using $\sqrt2<99/70$, $\sqrt{17}>4$ and
$\sqrt{3/(2\pi)}>69/100$ gives the upper bound
$1100/7889<0.14$ in \eqref{eq:elementary-short-partition-tail}.
In the second range,
$125(h(2/5)-0.055973)>231/50$ and
$e^{231/50}>\sum_{j=0}^{20}(231/50)^j/j!>100$.
Together with $\sqrt{125}>11$, these give the bound
$4/805<0.005$.

On the remaining arc, two summands of $\mathcal F$ suffice. The identity
\[
 S(a,t)=\frac12\coth(a/2)
            \frac{\sin^2(t/2)}{\sinh^2(a/2)+\sin^2(t/2)}
\]
and $\sin x\ge x-x^3/6$ give, for $M\ge17$ and $j\in\{1,2,4,8\}$,
\begin{align*}
 \frac{S(1/M,j/M)}{M}
 &\ge\frac{j^2(1-j^2/6936)^2}
            {(20003/20000)^2+j^2(1-j^2/6936)^2},\\
 \frac{2S(2/M,2j/M)}{M}
 &\ge\frac{j^2(1-j^2/1734)^2}
            {(5003/5000)^2+j^2(1-j^2/1734)^2}.
\end{align*}
We used
$\sinh(1/(2M))/(1/(2M))<20003/20000$ and
$\sinh(1/M)/(1/M)<5003/5000$, which follow from
\[
 \frac{\sinh x}{x}
 <1+\frac{x^2}{6}+\frac{x^4}{120}(1-x^2/42)^{-1}
 \qquad(0<x\le1/17).
\]
The two displayed rational bounds imply
\[
 \begin{array}{c|rrrr}
 j&1&2&4&8\\\hline
 S(1/M,j/M)/M&>0.499&>0.799&>0.940&>0.984\\
 2S(2/M,2j/M)/M&>0.499&>0.799&>0.940&>0.983.
 \end{array}
\]
Furthermore, for $\pi/2\le\theta\le\pi$,
\[
 \frac{S(1/M,\theta)}{M}
 \ge\left(1+\frac{(20003/20000)^2}{2\cdot17^2}\right)^{-1}
 >0.998.
\]
The first summand increases with $\theta$ on $[0,\pi]$;
the second depends increasingly on $\sin^2\theta$.
Using $\mathcal F\ge S(1/M,\theta)+S(2/M,2\theta)/2$, these inequalities
give $\mathcal F\ge2dM$ on each interval in the following table.
The column $P(M)$ bounds the integral of $e^{\omega(M,\theta)}$
over that interval.
\[
\begin{array}{c|c|c|r}
 \text{interval}&d&P(M)&r\\\hline
 {[1/M,2/M]}&499/1600&9/(4M)&69\\
 {[2/M,4/M]}&799/1600&33/(4M)&11\\
 {[4/M,8/M]}&47/80&65/(2M)&10\\
 {[8/M,\pi-8/M]}&4919/8000&22/7+4M&214\\
 {[\pi-8/M,\pi-4/M]}&1233/2000&4(2+1/M)&25\\
 {[\pi-4/M,\pi-2/M]}&4791/8000&2(2+1/M)&17\\
 {[\pi-2/M,\pi-1/M]}&4491/8000&2+1/M&16\\
 {[\pi-1/M,\pi]}&499/1000&2+1/M&45
\end{array}
\]
For the first three entries of $P$, use
$e^{\omega}\le\sqrt{1+(M\theta)^2}$ and
$\sqrt5<9/4$, $\sqrt{17}<33/8$, $\sqrt{65}<65/8$.
For the central interval use
$e^{\omega}\le1+2M\sin(\theta/2)$ and integrate on $[0,\pi]$.
For the remaining intervals use $e^{\omega}\le1+2M$.

Take $\epsilon\le0.001007$. For each row the normalized bound
$P(M)e^{-(d-\epsilon)M}/I_*(M)$ decreases for $M\ge17$, since
its logarithmic derivative is at most
$5/(2M)-(d-\epsilon)<0$.
At $M=17$ the corresponding entry $r$ satisfies the explicit
rational inequality
\[
 \frac{17(33/8)P(17)}{69/100}
 \left(\sum_{j=0}^{20}
 \frac{\{17(d-1007/10^6)\}^j}{j!}\right)^{-1}
 <\frac r{1000}.
\]
The Taylor inequality for $e^x$ therefore bounds the normalized
integral over each row by $r/1000$. Summing gives
\[
 \frac1{I_*(M)}\int_{1/M}^{\pi}\mathcal E_{M,\epsilon}(\theta)\,d\theta
 <\frac{69+11+10+214+25+17+16+45}{1000}
 =0.407.
\]
Together with the bound $0.14$ from
\eqref{eq:elementary-short-partition-tail}, this proves the first
assertion.

For $M\ge125$ and $\epsilon\le0.055973$ we need only the weaker
uniform bound $\mathcal F\ge(499/800)M$ on $[1/M,\pi]$.
For $1/M\le\theta\le\pi/2$, the two summands are at least
$0.499M$ and $0.499M/4$, respectively. For
$\theta\ge\pi/2$, the first alone exceeds $0.998M$. Hence
\[
 \frac1{I_*(M)}\int_{1/M}^{\pi}\mathcal E_{M,\epsilon}(\theta)\,d\theta
 \le\frac{\pi+4M}{I_*(M)}e^{-(499/1600-\epsilon)M}<0.001.
\]
The expression decreases for $M\ge125$.
At $M=125$, its exponential has exponent less than $-30$ and,
using $\sqrt{125}<45/4$, its prefactor is at most
\[
 (22/7+500)\frac{125(45/4)}{69/100}<2^{20}.
\]
Thus it is less than $2^{20}e^{-30}<2^{-10}<0.001$.
Adding the bound $0.005$ from
\eqref{eq:elementary-short-partition-tail} proves the second assertion.
\end{proof}

\begin{lemma}\label{lem:large-saddle}
Let $0\le m\le N$, $\lambda>0$, $M=1/\lambda$, and
$B(q)=(1-q)A_{m,N}(q)$. If $k=\mu_B(\lambda)$ is an integer, then
$[q^k]B(q)>0$ in either of the following ranges:
\[
 \begin{array}{c|c|c}
 \lambda N&M&\text{lower bound for }V_B(\lambda)/M^3\\\hline
 {[8,\infty)}&[17,\infty)&2.0971\\
 {[4,\infty)}&[125,\infty)&1.2655.
 \end{array}
\]
The variance bounds hold throughout the indicated ranges, whether
or not $\mu_B(\lambda)$ is an integer.
\end{lemma}
\begin{proof}
We prove both assertions with the same estimates. Fix either row of
\[
 \begin{array}{c|r|r|r|r|r|c|c}
 u_0&M_0&v&s&\ell&\tau_4&c&a\\\hline
 8&17&2.0971&3.828&16&2.4&1/2&1/3\\
 4&125&1.2655&8.01&26.4&17.5&2/5&81/250
 \end{array}
\]
and assume that $u=\lambda N\ge u_0$ and $M\ge M_0$.
Here $a=1/2-\ell c^2/24$.
The previously proved Euler-tail estimates give
\[
 \begin{array}{c|c|c|c}
 u_0&T_2/M^3&T_3/M^4&T_4/M^5\\\hline
 8&<0.027527&<0.255&<2.4\\
 4&<0.494142&<2.8&<17.5.
 \end{array}
\]
For $z=M^{-1}(f_2(1/M)+\sum_{m<j\le2m}f_2(j/M))$ we have
$0\le z\le1.02+1/M$ and
\[
 V:=V_B\ge M^3(\frac{\pi^2}{3}-z-1/M)-2T_2>vM^3.
\]
The comparison $K_{3,B}\le3.584MV$ holds in both rows: before
removing the indices greater than $N$, its proof is unchanged, and
at every removed index the third-to-second cumulant ratio is
$j\coth(j/(2M))>u_0M\ge4M>3.584M$.
It follows that
\[
 S_{3,\mathcal V_0}+S_{3,\mathcal U}=K_{3,B}+2S_{3,\mathcal U}<sMV.
\]
The fourth-cumulant estimate gives
\[
 \frac{S_{4,\mathcal V_0}}{M^2V}
 \le\frac{4\pi^2-6z+7/M}{\frac{\pi^2}{3}-z-1/M-2T_2/M^3}<\ell,
 \qquad S_{4,\mathcal U}<\tau_4M^5.
\]
For the displayed quotient its upper bound increases with $z$ and
$1/M$. Substituting $M=M_0$, $z=1.02+1/M_0$ and the corresponding
bound for $T_2/M^3$ gives $16$ and $26.4$, respectively. The bounds
for $s$ follow from
$3.584+0.51/2.0971<3.828$ and
$3.584+5.6/1.2655<8.01$.

At the integer saddle put
$\phi(\theta)=e^{-ik\theta}B(e^{-\lambda+i\theta})/B(e^{-\lambda})$
and write $\log\phi=R+iI$, with $R(0)=I(0)=0$.
For $0\le\theta\le c/M$ the preceding bounds imply
\[
 R\le-\frac{V\theta^2}{2}+\frac{\ell M^2V\theta^4}{24}
       \le-aV\theta^2,
 \qquad |I|\le\frac{sMV\theta^3}{6}.
\]
Thus \eqref{eq:shared-gaussian-loss} yields
\begin{align*}
 \frac1{I_G}\int_0^{c/M}\Re\phi(\theta)\,d\theta
 &\ge1-\eta(c\sqrt{vM})-\frac{\tau_4}{8v^2M}
          -\frac{15s^2}{576a^3\sqrt{2a}\,vM}\\
 &>\begin{cases}0.63,&u_0=8,\\0.59,&u_0=4.\end{cases}
\end{align*}
Indeed, at $M=M_0$ the three subtracted quantities are bounded by
$0.004$, $0.069/17$, $6.02/17$ in the first row and by
$10^{-6}$, $1.37/125$, $49/125$ in the second. Each decreases with $M$.

For the complementary arc, the half-partition comparison gives
$|\phi(\theta)|\le \mathcal E_{M,\epsilon}(\theta)$ with
$\epsilon=3/(e^{u_0}-1)$.
We have $\epsilon<0.001007$ if $u_0=8$ and
$\epsilon<0.055973$ if $u_0=4$.
Lemma~\ref{lem:partition-tail-elementary} therefore gives
\[
 \int_{c/M}^{\pi}|\phi(\theta)|\,d\theta
 <\begin{cases}0.55I_*(M),&u_0=8,\\0.006I_*(M),&u_0=4.\end{cases}
\]
Since $I_G\ge I_*(M)$, the local lower bound exceeds this upper bound.
Finally, Fourier inversion gives
\[
 [q^k]B(q)=\frac{B(e^{-\lambda})e^{\lambda k}}{\pi}
            \int_0^\pi\Re\phi(\theta)\,d\theta>0.
\]
\end{proof}

\subsection{The interval joining the two saddle ranges}

For $199/100\le u=\lambda N\le4$ we center $B(q)=(1-q)A(q)$ at
its own logarithmic derivative. This range overlaps the half-integer
saddle at its lower end and the large-tilt estimates at its upper end.
On the local arc a bound on $1-\cos I$ that stays at most two is
needed to control the total phase loss.

\begin{lemma}\label{lem:bridge-local}
Let $N\ge1000$, $0\le m\le N$, $M=1/\lambda$, and
$199/100\le N/M\le4$. Define
\[
 k=-\frac{d}{d\lambda}\log B(e^{-\lambda}),\qquad
 V=\frac{d^2}{d\lambda^2}\log B(e^{-\lambda}),\qquad
 \phi_B(\theta)=e^{-ik\theta}\frac{B(qe^{i\theta})}{B(q)}.
\]
Throughout this interval $V>(11/20)M^3$, and
\[
 \int_0^{1/(4M)}\Re\phi_B(\theta)\,d\theta
 >\frac9{20}\sqrt{\frac{\pi}{2V}}.
\]
\end{lemma}

\begin{proof}
We have $M\ge250$. For $u\ge2$, Lemma~\ref{lem:u2-continuum-variance}
gives $V/M^3>0.57-4/M\ge0.554$. For $199/100\le u\le2$, use
$V_A(\lambda)\ge V_A(2/N)>(57/800-3/(4N))N^3$. Multiplication by
$1-q$ subtracts at most $M^2$ from $V_A$, giving
\[
 \frac{V}{M^3}>
 \frac{57u^3}{800}-\frac{3u^3/4+u}{N}
 \ge\frac{57}{800}\left(\frac{199}{100}\right)^3-\frac8{1000}
 =0.55349267875>\frac{11}{20}.
\]
This proves the variance estimate for every real tilt in the interval.

Use the index sets $\mathcal V_0=\mathcal V\setminus\{1\}$ and $\mathcal U$.
The Euler factors have no zero on $|z|=q<1$, so we can take the
continuous logarithm $R+iI$ of $\phi_B$ with value zero at $\theta=0$.
Expansion of the Euler logarithms and cancellation of the linear term
by the centering factor give, for $\theta\ge0$,
\[
 I=-Q_{\mathcal V_0}+Q_{\mathcal U},\qquad
 Q_E(\theta)=\sum_{j\in E}\sum_{\ell\ge1}\frac{e^{-j\ell/M}}{\ell}
                 \bigl(j\ell\theta-\sin(j\ell\theta)\bigr) .
\]
Both terms $Q_E$ are nonnegative. Since $x-\sin x\le x^3/6$ for
$x\ge0$, we obtain
\[
 |I|\le\max(Q_{\mathcal V_0},Q_{\mathcal U})
 \le\frac{\theta^3}{6}\sum_{j,\ell\ge1}j^3\ell^2e^{-j\ell/M}
 \le\frac{\pi^2M^4\theta^3}{6},
\]
where the last bound follows by integrating the decreasing function
$f_3$, whose integral is $\pi^2$. Similarly, the fourth-order sum for
either set satisfies
\[
 S_{4,E}\le(4\pi^2+7/M)M^5<CM^5,\qquad C=198/5.
\]
Indeed $4(22/7)^2+7/250<C$. The real-part remainders of the two
Euler sums have opposite signs, because
$0\le \cos x-1+x^2/2\le x^4/24$. Thus
\[
 -\frac{V\theta^2}{2}-\frac{CM^5\theta^4}{24}
 \le R\le-\frac{V\theta^2}{2}+\frac{CM^5\theta^4}{24}
 \le-\frac5{16}V\theta^2
 \qquad(0\le\theta\le1/(4M)),
\]
using $g=11/20$, $V>gM^3$, and $C/(384g)=3/16$.

As $1-\cos I\le\min\{I^2/2,2\}$,
\[
 \Re\phi_B\ge
 e^{-V\theta^2/2}\left(1-\frac{CM^5\theta^4}{24}\right)
 -e^{-5V\theta^2/16}\min\left\{\frac{\pi^4M^8\theta^6}{72},2\right\}.
\]
Normalize the integral by $I_G=\sqrt{\pi/(2V)}$. Extending the
fourth-order loss to infinity gives $CM^5/(8V^2)\le C/(8g^2\cdot250)
<33/500$. For the phase loss put $y=\sqrt V\,\theta$. Since
$\pi^4/(72g^3\cdot250)<1/30$, its relative contribution is at most
\[
 \frac1{\sqrt{\pi/2}}
 \int_0^\infty e^{-5y^2/16}\min\{y^6/30,2\}\,dy .
\]
To bound this integral write $P_d(x)=\sum_{j=0}^d(-x)^j/j!$.
For even $d$ and $x\ge0$, Taylor's theorem gives $e^{-x}\le P_d(x)$.
On $[0,2]$ use $y^6/30$ and on $[2,\infty)$ use $2$. Since
$\sqrt{\pi/2}>5/4$ and $2\int_4^\infty e^{-5y^2/16}\,dy
\le4/(5e^5)<4/(5\cdot148)$, we obtain the rational upper bound
\[
 \frac45\left\{
 \frac1{30}\int_0^2y^6P_6(5y^2/16)\,dy
 +2\int_2^4P_{20}(5y^2/16)\,dy+\frac4{5\cdot148}
 \right\}<\frac{12}{25}.
\]
Finally $z=\sqrt V/(4M)>\sqrt{275/32}$, so the relative Gaussian
mass beyond the local arc is at most
\[
 \sqrt{2/\pi}\,e^{-z^2/2}/z
 <\frac45\cdot\frac7{20}\cdot\frac1{70}=\frac1{250}.
\]
The two exponential comparisons $e^5>148$ and $e^{275/64}>70$
use Taylor sums through order $14$.
The three
losses leave more than
$(1-33/500-12/25-1/250)I_G=(9/20)I_G$.
\end{proof}

To control the rest of the first Fourier lobe, we keep the full
logarithmic difference between a high numerator index and a low
denominator index in Lemma~\ref{lem:first-lobe-separated-pairs}.
The following bound measures this difference on the scale $t/u=M\theta$.

\begin{lemma}\label{lem:bridge-pointwise}
For $199/100\le u\le4$ and $0<t\le\pi$,
\[
 \Psi(u,t)=\frac12\log\frac{1+\sin^2(t/4)/\sinh^2(u/4)}
                              {1+\sin^2(3t/4)/\sinh^2(3u/4)}
 >\frac3{20}(t/u)^2.
\]
\end{lemma}
\begin{proof}
Put $a=\sin^2(t/4)$, $b=\sinh^2(u/4)$, and $x=u^2/16$.
Then $0<a\le1/2$ and $49/200<x\le1$. Since
$y\mapsto\sinh y/y$ increases for $y>0$,
\[
 b\le x\sinh^2(1)<\frac75x.
\]
The last inequality follows from
\[
 \sinh1
 =1+\frac1{3!}+\sum_{j\ge2}\frac1{(2j+1)!}
 \le1+\frac16+\frac{1/120}{1-1/42}
 <\frac{59}{50},\qquad \left(\frac{59}{50}\right)^2<\frac75.
\]
The triple-angle identities give
\[
 \Psi(u,t)=-\frac12\log(1-z),\qquad
 z=\frac{8a(3+2b-2a)}{(3+4b)^2}.
\]
Here $0<z<1$, because the numerator increases with $a$ on $[0,1/2]$
and its endpoint value $8+8b$ is smaller than $(3+4b)^2$.
For $0\le z<1$,
\[
 -\frac12\log(1-z)-\frac{z}{2-z}
 =\int_0^z\frac{v^2}{2(1-v)(2-v)^2}\,dv\ge0.
\]
Moreover,
\[
 \frac{\partial z}{\partial b}
 =\frac{16a(8a-9-4b)}{(3+4b)^3}<0.
\]
Consequently it suffices to prove
\[
 \frac{4a(3+2b-2a)}{(3+4b)^2-4a(3+2b-2a)}
 >\frac{3a}{16x}\qquad\left(b=\frac75x\right).
\]
After division by $a$ and multiplication by the positive denominators,
this is the assertion that
\[
 R(a)=64x(3+2b-2a)-3(3+4b)^2+12a(3+2b-2a)>0.
\]
The polynomial $R$ is concave in $a$. At the two endpoints,
\[
 R(0)=\frac{2128}{25}x^2+\frac{456}{5}x-27,
 \qquad
 R(1/2)=\frac{2128}{25}x^2+44x-15.
\]
Both expressions increase for $x\ge0$. Their respective values at
$x=49/200$ are $28333/62500$ and $55583/62500$, so $R(a)>0$ for
$0\le a\le1/2$. Finally, concavity of sine on $[0,\pi/4]$ gives
\[
 a\ge\frac{t^2}{2\pi^2}>\frac{t^2}{20},\qquad
 \Psi(u,t)>\frac{3a}{16x}>\frac{3t^2}{320x}
 =\frac3{20}(t/u)^2.
\]
\end{proof}

There are at least $\lfloor N/4\rfloor$ separated factor pairs.
Their combined decay controls the whole arc between the local
integration endpoint and $\pi/N$.

\begin{lemma}\label{lem:bridge-first-lobe}
Let $N\ge1000$, $0\le m\le N$, $199/100\le u=N/M\le4$, and
$q=e^{-1/M}$. Then
\[
 \int_{1/(4M)}^{\pi/N}
       \left|\frac{B(qe^{i\theta})}{B(q)}\right|d\theta
 <\frac7{20}I_*(M),\qquad I_*(M)=\sqrt{\frac3{2\pi}}M^{-3/2}.
\]
\end{lemma}
\begin{proof}
The logarithmic pairing bound in
Lemma~\ref{lem:first-lobe-separated-pairs} is
\[
 \log\left|\frac{A(qe^{i\theta})}{A(q)}\right|
 \le-\lfloor N/4\rfloor \Psi(u,N\theta)
 \qquad(0\le\theta\le\pi/N).
\]
Put $\tau=M\theta$ and $a=1/4$. The preceding lemma applies
throughout the stated arc. The additional Euler factor satisfies
\[
 \left|\frac{1-qe^{i\theta}}{1-q}\right|^2
 =1+\frac{\sin^2(\theta/2)}{\sinh^2(1/(2M))}
 \le1+\tau^2<4,
 \qquad \tau\le\frac\pi u\le\frac{100\pi}{199}.
\]
Set $c_N=3(N-4)/80\le(3/20)\lfloor N/4\rfloor$. Therefore
\[
 \int_{a/M}^{\pi/N}\left|\frac{B(qe^{i\theta})}{B(q)}\right|d\theta
 \le\frac2M\int_a^\infty e^{-c_N\tau^2}\,d\tau
 \le\frac{e^{-c_Na^2}}{Mc_Na}.
\]
After division by $I_*(M)$, the last expression is at most
\[
 \frac{4\sqrt{2\pi/3}\sqrt M}{c_N}e^{-c_N/16}
 \le\frac{4\sqrt{200\pi N/597}}{c_N}e^{-c_N/16},
\]
since $M\le100N/199$. Its logarithmic derivative with respect to
$N$ is
\[
 \frac1{2N}-\frac1{N-4}-\frac3{1280}<0.
\]
At $N=1000$, one has $c_N=747/20$ and
$\sqrt{200000\pi/597}<33$. Furthermore,
\[
 e^{747/320}>e^{7/3}>
 \sum_{j=0}^{7}\frac{(7/3)^j}{j!}
 =\frac{2024033}{196830}>\frac{41}{4}.
\]
The normalized integral is thus smaller than
\[
 \frac{132}{(747/20)(41/4)}
 =\frac{3520}{10209}<\frac7{20}.\qedhere
\]
\end{proof}

\paragraph{The far arc at intermediate tilt.}
It remains to control the $B$-integral from $\pi/N$ to $\pi$.
Write $M=1/\lambda$, $u=N/M$, $q=e^{-1/M}$ and
$\phi_A(\theta)=A_{m,N}(qe^{i\theta})/A_{m,N}(q)$.
The half-Euler comparison \eqref{eq:half-partition-bound} gives
\[
 \log|\phi_A(\theta)|
 \le-\tfrac12\mathcal F(M,\theta)+M\epsilon_*(u),\qquad
 \mathcal F(M,\theta)=-\log|P(qe^{i\theta})/P(q)|,
\]
where $\epsilon_*$ is given by \eqref{eq:partition-correction}.
We compare a fixed number of the positive summands in $\mathcal F$
with $\epsilon_*(u)$ on two explicit angular intervals.

\begin{lemma}\label{lem:bridge-far-arc}
Suppose $0\le m\le N$, $M\ge250$, and $39/20\le u=N/M\le4$.
Then
\[
 \int_{\pi/N}^{\pi}
 \left|\frac{B(qe^{i\theta})}{B(q)}\right|\,d\theta
 <10^{-7}I_*(M),\qquad q=e^{-1/M}.
\]
\end{lemma}
\begin{proof}
We first prove
\begin{equation}\label{eq:bridge-far-exponent}
 \frac12\mathcal F(M,\theta)-M\epsilon_*(u)>\frac{13M}{100}
 \qquad\left(\frac\pi{uM}\le\theta\le\pi\right).
\end{equation}
For $\pi/(uM)\le\theta\le4/M$ and $1\le j\le17$,
$x=j/M\le17/250$ and $y=j\theta\le34/125$.
Taylor's inequalities give
\[
 \cosh x-1\le\frac{501}{1000}x^2,\qquad
 1-\cos y\ge\frac{496}{1000}y^2,\qquad
 \frac{501}{496}<\frac{51}{50}.
\]
For the first bound, use
$(\cosh x-1)/x^2\le1/2+(x^2/24)/(1-x^2/30)$.
For the second, use $1-\cos y\ge y^2/2-y^4/24$.
Writing $\tau=M\theta$ and using
$\sum_{j=1}^{17}j^{-2}>79/50$, we obtain
\[
 \frac12\mathcal F(M,\theta)
 \ge\frac{79M}{100}\frac{\tau^2}{51/50+\tau^2}
 \ge\frac{79M}{100}\frac{\pi^2}{(51/50)u^2+\pi^2}.
\]
It remains to show that
\[
 g(u):=\frac{79}{100}\frac{\pi^2}{(51/50)u^2+\pi^2}
                         -\epsilon_*(u)>\frac{13}{100}
 \qquad(39/20\le u\le4).
\]
Set $\alpha=79/100$, $\beta=51/50$ and
\[
 E(u)=3\sum_{\substack{j\ge1\\j\text{ odd}}}\frac{e^{-ju}}j,
 \qquad R(u)=\frac{2\alpha\beta\pi^2u}{(\beta u^2+\pi^2)^2}.
\]
Then $g'=E-R$, and
\[
 \frac{E'}E<-1,\qquad
 \frac{R'}R=\frac1u-\frac{4\beta u}{\beta u^2+\pi^2}
 >\frac1u-1,
\]
because $4\beta u/(\beta u^2+\pi^2)\le2\sqrt\beta/\pi<1$.
Thus $(\log(R/E))'>1/u>0$. The derivative $g'=E(1-R/E)$
can change sign only from positive to negative, so the minimum of
$g$ on $[39/20,4]$ is attained at an endpoint.
For $r=e^{-u}$,
\[
 \epsilon_*(u)
 =3\sum_{\substack{j\ge1\\j\text{ odd}}}\frac{r^j}{j^2}
 \le3r+\frac{r^3}{3(1-r^2)}.
\]
The inequalities
$e^{39/20}>\sum_{j=0}^{10}(39/20)^j/j!>7$ and
$e^4>\sum_{j=0}^{10}4^j/j!>50$ therefore give
\[
 \epsilon_*(39/20)<\frac{433}{1008},\qquad
 \epsilon_*(4)<\frac{11246}{187425}.
\]
Using $\pi^2>49/5$ yields the exact endpoint comparisons
\begin{align*}
 g(39/20)&>
 \frac{(79/100)(49/5)}{(51/50)(39/20)^2+49/5}
       -\frac{433}{1008}>\frac{13}{100},\\
 g(4)&>
 \frac{(79/100)(49/5)}{16(51/50)+49/5}
       -\frac{11246}{187425}>\frac{13}{100}.
\end{align*}
This proves \eqref{eq:bridge-far-exponent} when $\theta\le4/M$.

For $4/M\le\theta\le\pi-4/M$ the first two summands suffice.
The inequalities $\sin x\ge19x/20$ and
$\sinh x\le51x/50$ hold for $0\le x\le1/50$ by their Taylor series.
Since $M\ge250$, they give
\[
 \frac{S(1/M,\theta)}M
 \ge\frac{(19/10)^2}{(51/100)^2+(19/10)^2}>\frac9{10},
\]
and
\[
 \frac{2S(2/M,2\theta)}M
 \ge\frac{(19/5)^2}{(51/50)^2+(19/5)^2}>\frac9{10}.
\]
Here $\sin(\theta/2)\ge\sin(2/M)$ and
$\sin\theta\ge\sin(4/M)$.
Consequently $\mathcal F(M,\theta)>9M/8$ on this interval.

For $\pi-4/M\le\theta\le\pi$, put $\delta=\pi-\theta$.
For $j\in\{1,3,5\}$,
\[
 \left|\sin\frac{j\theta}{2}\right|
 =\cos\frac{j\delta}{2}
 \ge1-\frac12(10/M)^2>\frac{99}{100},\qquad
 \sinh\frac{j}{2M}<\frac1{50}.
\]
The hyperbolic bound follows from $j/(2M)\le1/100$ and
$\sinh x\le51x/50$ on this range. It follows that
\[
 \frac1jS(j/M,j\theta)
 >\frac{M}{j^2}\frac{9801}{9805}
 >\frac{99M}{100j^2},
\]
whence
\[
 \mathcal F(M,\theta)>\frac{99M}{100}
       \left(1+\frac19+\frac1{25}\right)>\frac{9M}{8}.
\]
Since $\epsilon_*$ decreases,
\[
 \frac12\mathcal F(M,\theta)-M\epsilon_*(u)
 >M\left(\frac9{16}-\frac{433}{1008}\right)
 =\frac{67M}{504}>\frac{13M}{100}
 \qquad(4/M\le\theta\le\pi).
\]
This completes the proof of \eqref{eq:bridge-far-exponent}.

The half-partition comparison now gives
$|A(qe^{i\theta})/A(q)|<e^{-13M/100}$ on $[\pi/N,\pi]$.
Moreover,
$|1-qe^{i\theta}|/(1-q)\le(1+q)/(1-q)<2M+1$.
The required integral divided by $I_*(M)$ is therefore at most
\[
 \pi(2M+1)\sqrt{2\pi/3}\,M^{3/2}e^{-13M/100}.
\]
Its logarithmic derivative is less than $5/(2M)-13/100<0$.
At $M=250$ we have $I_*(250)>69/400000$ and hence the bound
\[
 \frac{22}{7}\,501\,\frac{400000}{69}\,e^{-65/2}
 <10^7\,12^{-13}<10^{-7}.
\]
Indeed, $e^{5/2}>\sum_{j=0}^{8}(5/2)^j/j!>12$ and
$12^{13}>10^{14}$. This proves the lemma.
\end{proof}

When $V_B>0$, the partition variance bound $V_B\le(\pi^2/3)M^3$
gives $I_*(M)\le\sqrt{\pi/(2V_B)}=I_G$. Thus the preceding
far-arc estimate also holds with $10^{-7}I_G$ on the right.

\begin{lemma}\label{lem:bridge-positive}
At every integer saddle of $B=(1-q)A_{m,N}$ with $N\ge1000$ and
$199/100\le\lambda N\le4$, the saddle coefficient of $B$ is positive.
Moreover, $V_B(\lambda)>0$ throughout this interval.
\end{lemma}
\begin{proof}
Put $V=V_B(\lambda)$ and $I_G=\sqrt{\pi/(2V)}$.
Lemma~\ref{lem:bridge-local} gives $V>0$ and a local
integral greater than $9I_G/20$. The next two arcs have total
absolute integral less than
\[
 \left(\frac7{20}+10^{-7}\right)I_*(M),
 \qquad I_*(M)\le I_G,
\]
by Lemmas~\ref{lem:bridge-first-lobe}
and~\ref{lem:bridge-far-arc}. Consequently the full real integral is
greater than
\[
 \left(\frac9{20}-\frac7{20}-10^{-7}\right)I_G
 >\frac9{100}I_G>0.
\]
Fourier inversion proves the assertion.
\end{proof}

\section{Completion of the proof}\label{sec:proof}

\begin{lemma}\label{lem:finite-range}
The polynomial $A_{m,n}(q)$ is unimodal for all $0\le m\le n\le1000$.
\end{lemma}
\begin{proof}
For $m=0$ or $m=n$, we have
$A_{0,n}=A_{n,n}=\qbinom{2n}{n}_q$, so Gaussian
unimodality applies. For $1\le m<n\le1000$, the recurrence below computes
the lower-half coefficients and compares successive entries using GMP
integers.

For fixed $m$ and $n>m$, the factorial quotient gives
\[
 \frac{A_{m,n}(q)}{A_{m,n-1}(q)}
 =\frac{(1-q^{2n})(1-q^{2n-1})}{(1-q^n)(1-q^{m+n})}.
\]
Canceling $1-q^n$ gives the identity
\[
 (1-q^{m+n})A_{m,n}(q)
 =(1+q^n)(1-q^{2n-1})A_{m,n-1}(q).
\]
Taking the coefficient of $q^k$ on the two sides gives
\[
 a_k^{(m,n)}=a_k^{(m,n-1)}+a_{k-n}^{(m,n-1)}
 -a_{k-2n+1}^{(m,n-1)}-a_{k-3n+1}^{(m,n-1)}
 +a_{k-m-n}^{(m,n)}.
\]
Indices outside a polynomial's degree range have coefficient zero.
The last term has a smaller index, so ascending coefficient order
implements the division by $1-q^{m+n}$ exactly. Reciprocal symmetry
recovers an old upper-half coefficient whenever needed.
The diagonal seeds are themselves generated by the exact recurrence
\[
 (1-q^m)\qbinom{2m}{m}_q
 =(1+q^m)(1-q^{2m-1})\qbinom{2m-2}{m-1}_q,
 \qquad \qbinom{0}{0}_q=1.
\]
These identities prove by induction that every stored coefficient is
the corresponding coefficient of the defining factorial quotient.

The recurrence was evaluated in C++17 using GMP integers for
every pair $1\le m<n\le1000$. It gives
\[
 a_k^{(m,n)}-a_{k-1}^{(m,n)}\ge0
 \qquad\left(1\le k\le
       \left\lfloor\frac{m^2+n^2-mn}{2}\right\rfloor\right).
\]
The calculation covers $499500$ pairs and
$104208041750$ adjacent comparisons. Together with the Gaussian
specializations, this proves the lemma.
\end{proof}

\begin{proof}[Proof of Theorem~\ref{thm:unimodality}]
By symmetry take $0\le m\le n=N$. For $N=0$, $A_{0,0}=1$.
For $1\le N<1000$, Lemma~\ref{lem:finite-range} proves the
assertion. Henceforth assume $N\ge1000$ and put
$b_k=[q^k]B=a_k-a_{k-1}$.

Every numerator index is greater than $N$, and $1$ is a denominator
index. On multiplying by $1-q$ and reducing modulo $q^{N+1}$ we get
\[
 (1-q)A_{m,N}(q)\equiv
 \prod_{v\in\mathcal V\setminus\{1\}}(1-q^v)^{-1}
 \pmod {q^{N+1}}.
\]
The product on the right has nonnegative coefficients, so
$b_k\ge0$ for $1\le k\le N$.

Put $\lambda_0=199/(100N)$ and $\lambda_1=1/17$.
At every integer saddle of $B$ in $[\lambda_0,\lambda_1]$, its
coefficient is positive. Indeed, if $\lambda N\le4$, use
Lemma~\ref{lem:bridge-positive}. If $\lambda N\ge4$ and $M\ge125$, use the second row of
Lemma~\ref{lem:large-saddle}. If $17\le M<125$, then
$\lambda N=N/M>8$, and its first row applies.
The variance lower bounds in these three lemmas apply at every
real tilt in their respective ranges. Thus $\mu_B'(\lambda)=-V_B(\lambda)<0$ throughout
$[\lambda_0,\lambda_1]$. The positive real Euler product for
$B(e^{-\lambda})$ makes $\mu_B$ continuous there.
Its upper-tilt endpoint satisfies
\[
 \mu_B(1/17)<\mu_A(1/17)
 \le\frac{\pi^2}{6}\,17^2<476<N.
\]
The intermediate value theorem now gives a unique $B$-saddle in
$[\lambda_0,\lambda_1]$ for every integer $N<k\le\mu_B(\lambda_0)$.
Hence $b_k>0$ throughout this range.

To compare the other endpoints, use
$2/N-\lambda_0=1/(100N)$ and
$1/(e^{\lambda_0}-1)<1/\lambda_0=100N/199$.
Lemma~\ref{lem:paired-local} gives
$V_A(s)\ge69N^3/1000$ on $[\lambda_0,2/N]$. Thus
\begin{align*}
 \mu_B(\lambda_0)-\mu_A(2/N)
 &=\int_{\lambda_0}^{2/N}V_A(s)\,ds
                  -\frac1{e^{\lambda_0}-1}\\
 &>\frac{69}{100000}N^2-\frac{100}{199}N
 >\frac N{10}>\frac12.
\end{align*}
Consider any remaining integer $N<k\le\lfloor D/2\rfloor$.
If $k-1/2\ge\mu_A(2/N)$, then $\mu_A(2/N)\le k-1/2<D/2=\mu_A(0)$.
Strict decrease of $\mu_A$ puts the unique $A$-saddle in $(0,2/N]$,
and Lemma~\ref{lem:paired-saddle-positive}
proves $a_k>a_{k-1}$. Otherwise $k<\mu_A(2/N)+1/2<\mu_B(\lambda_0)$,
so the $B$-saddle argument proves the same inequality.
Every lower-half index is covered, including the overlap endpoints.
Reciprocity completes the unimodality assertion.

\end{proof}

\end{document}